\documentclass{amsart}
\usepackage{amsmath,amssymb,amsfonts}
\usepackage{amsthm}
\usepackage{url,graphicx,enumerate,color,hyperref,booktabs,epstopdf}
\usepackage[top=1in,bottom=1in,left=1in,right=1in]{geometry}
\newtheorem{definition}{Definition}[section]
\newtheorem{remark}{Remark}[section]

\numberwithin{equation}{section}
\definecolor{orange}{rgb}{1,0.5,0}
\definecolor{dgreen}{rgb}{0,0.5,0}
\definecolor{lightblue}{rgb}{0.27,0.51,0.706}
\newcommand{\Srm}{\textrm{\S}}

\newcommand{\dt}{\Delta t}
\newcommand{\dx}{\Delta x}
\newcommand{\myvec}[1]{\boldsymbol{#1}}
\newcommand{\inprod}[2]{\langle #1,#2\rangle} % define inner product <f,g>
\newcommand{\CURL}{\nabla\times}
\newcommand{\DIV}{\nabla\cdot}
\newcommand{\tdiv}{\textnormal{div}}
\newcommand{\tcurl}{\textnormal{curl}}
\newcommand{\intd}{\ \mathrm{d}}
\newcommand{\nl}[3]{(\myvec{#1}^{#2}_{#3}\cdot\nabla)\myvec{#1}^{#2}_{#3}}
\newcommand{\NL}[2]{\myvec{N}(\myvec{#1}_{#2})}
\newcommand{\bdryint}[1]{\int_{\partial\Omega}{#1}\intd S}

\newcommand{\imex}{IMEX}

\usepackage{soul}

\begin{document}

\title[High-Order Methods for PPE Reformulation for NSE with EBC]
{High-order Methods for a Pressure Poisson Equation Reformulation of the
Navier-Stokes Equations\\ with Electric Boundary Conditions}

\author[R. R. Rosales]{Rodolfo Ruben Rosales}
\address{Department of Mathematics,
Massachusetts Institute of Technology,
77 Massachusetts Avenue, Cambridge, MA 02139}
\email{rrr@math.mit.edu}

\author[B. Seibold]{Benjamin Seibold}
\address{Department of Mathematics, Temple University,
1805 North Broad Street, Philadelphia, PA 19122}
\email{seibold@temple.edu}
\urladdr{http://www.math.temple.edu/~seibold}

\author[D. Shirokoff]{David Shirokoff}
\address{Department of Mathematical Sciences,
New Jersey Institute of Technology,
University Heights, Newark, NJ 07102}
\email{david.g.shirokoff@njit.edu}
\urladdr{http://web.njit.edu/~shirokof}

\author[D. Zhou]{Dong Zhou}
\address{Department of Mathematics, California State University, Los Angeles,
5151 State University Drive, Los Angeles, CA 90032}
\email{dong.zhou@calstatela.edu}
%\urladdr{http://www.math.temple.edu/~dzhou}

\date{\today}
%\dedicatory{}
%----------------------------------------------------------------------------
% Abstract
%----------------------------------------------------------------------------
\begin{abstract}
Pressure Poisson equation (PPE) reformulations of the incompressible Navier-Stokes equations (NSE) replace the incompressibility constraint by a Poisson equation for the pressure and a suitable choice of boundary conditions. This yields a time-evolution equation for the velocity field only, with the pressure gradient acting as a nonlocal operator. Thus, numerical methods based on PPE reformulations, in principle, have no limitations in achieving high order. In this paper, it is studied to what extent high-order methods for the NSE can be obtained from a specific PPE reformulation with electric boundary conditions (EBC). To that end, implicit-explicit (\imex{}) time-stepping is used to decouple the pressure solve from the velocity update, while avoiding a parabolic time-step restriction; and mixed finite elements are used in space, to capture the structure imposed by the EBC. Via numerical examples, it is demonstrated that the methodology can yield at least third order accuracy in space and time.
\end{abstract}
%----------------------------------------------------------------------------
\maketitle
% --------------------------------------------------------------------------- 80
%\input{ppe_input}
% --------------------------------------------------------------------------- 80
\section{Introduction}
% --------------------------------------------------------------------------- 80
Developing efficient high-order time stepping methods for the incompressible
Navier-Stokes equations (NSE) is challenging due to the fact that the velocity
and pressure are coupled via an incompressibility constraint.
Numerical methods that treat both the velocity $\myvec{u}$ and the pressure
$p$ in a fully implicit fashion \cite{KrzywickiLadyzhenskaya1966} provide a
comparatively straightforward pathway towards high-order. However, fully implicit approaches result in large (possibly nonlinear) saddle-point
systems. This imposes a requirement to select spatial discretizations that ensure stable solutions of the resulting discrete equations (i.e.\, staggered
grids~\cite{HarlowWelch1965}, or finite elements satisfying the inf-sup
condition~\cite{FortinBrezzi1991}). In addition, the resulting discrete saddle point problems are non-trivial to solve efficiently~\cite{BenziGolubLiesen2005,GreifReesSzyld2017,BrennerLiSung2014}. 

On the other hand, many numerical approaches that decouple $\myvec{u}$ and $p$
result in smaller systems of equations with fewer coupled variables, and often times avoid a saddle-point structure. Hence, such methodologies are attractive for certain large-scale problems. Unfortunately, it is challenging to achieve high-order (in time) when decoupling velocity and pressure. This paper develops finite element approaches for certain pressure Poisson equation (PPE) reformulations of the Navier-Stokes equations that allow for a systematic pathway towards high order while decoupling velocity and pressure.

Numerical methods for the NSE that decouple velocity and pressure date
back to the late 1960s with the introduction of projection
methods~\cite{Chorin1968,Temam1969}.
The idea of projection methods is to first evolve the velocity without the
pressure, and then project the velocity back into the space of divergence-free
fields via a Poisson problem. Projection methods are one example of a larger class of methods known as fractional step methods \cite{Perot1993}.
Efforts have been made in the past few decades to improve the accuracy of
projection and fractional step methods to second order (in time), and higher,
\cite{BellColellaGlaz1989,KarniadakisIsraeliOrszag1991,%
KimMoin1985,BrownCortezMinion2001}.
In particular, accurate methods that go beyond second order are an
ongoing area of research \cite{RebholzXiao2017,LiuLiuPego2010}. 
In addition to the difficulties of achieving high order in time accuracy, projection methods carry the risk of producing numerical boundary layers via the Poisson equation for the pressure, which causes a degradation in spatial error convergence. An extensive overview of projection methods is given by Guermond, Minev, and Shen in \cite{GuermondMinevShen2006}, where projection methods and various improvements, as well as their theoretical and numerical convergence results are discussed.
Recent progress on generalizing the artificial compressibility method \cite{GuermondMinev2015,GuermondMinev2017} (which has traditionally been first order in time) has lead to alternative avenues that obtain high-order (beyond second order) in time schemes.

PPE reformulations of the Navier-Stokes equations \cite{%
GreshoSani1987,%
Henshaw1994,%
HenshawKreissReyna1994,%
HenshawPetersson2003,%
JohnstonLiu2002,%
JohnstonLiu2004,%
KleiserSchumann1980,%
Rempfer2006,%
SaniShenPironneauGresho2006,%
ShirokoffRosales2010, Zhang2016}
provide an alternate route towards devising
high-order in time numerical methods.  The basic idea underlying PPE
approaches is not to discretize the NSE directly, but instead to
(i)~reformulate the NSE into a system of PDEs with a Poisson equation for
    the pressure in lieu of the divergence constraint (the PPE system),
(ii)~devise boundary conditions that ensure that the new set of equations
     guarantees incompressibility~\cite{GreshoSani1987,%
KarniadakisIsraeliOrszag1991,Henshaw1994,SaniShenPironneauGresho2006},
and then (iii)~discretize and solve the resulting PPE system.
The solution to the Poisson equation for the pressure is designed
so that the PPE reformulation is equivalent to
the original NSE for solutions that are sufficiently smooth. In other words, at the continuous level,
solving the PPE is equivalent to solving the Navier-Stokes equations.
In a sequence of work by Henshaw et al.~\cite{Henshaw1994,HenshawKreiss1995,%
HenshawKreissReyna1994,HenshawPetersson2003},
the recovery of the pressure through the solution of a pressure Poisson
problem was done in a discrete setting (with finite differences), using the
boundary condition $\DIV\myvec{u}=0$
together with other numerical boundary conditions.
PPE formulations as PDEs were introduced in
\cite{JohnstonLiu2004} by Johnston and Liu, and in \cite{ShirokoffRosales2010}
by Shirokoff and Rosales. Once the PPE has been formulated as a set of
continuous PDEs, one may then examine a variety of different numerical
discretizations (in both space and time) to solve the PPE system (and hence,
equivalently, the NSE).

One important advantage of PPE reformulations compared to the standard form
of the NSE is that by solving the Poisson equation, the pressure can be viewed
as a global function of the velocity field $p = P(\myvec{u})$. The fact that
the pressure can be written as a function of the velocity enables
straight-forward numerical approaches for decoupling the velocity and
pressure, i.e.\, to enable implicit-explicit time-stepping strategies.
PPE reformulations then provide the possibility to devise high-order in
time numerical methods within a systematic framework.
Note also that (unlike the NSE) PPE reformulations are defined even if the
initial conditions are not incompressible
\cite{ShirokoffRosales2010,IgnatovaIyerKelliherPegoZarnescu2013},
which can be an important advantage when dealing with real data.
Drawbacks of PPE reformulations are:
(i)~For the numerical solution the (numerical) divergence field is not
    exactly zero.
(ii)~The pressure boundary conditions are typically complicated. Thus their
   interaction with the velocity is not straightforward to understand and
   analyze.

Here we focus on numerical methods based on the Shirokoff-Rosales (SR) PPE
reformulation proposed in \cite{ShirokoffRosales2010}.
Numerical discretizations for the SR PPE reformulation have been proposed
in the context of standard finite difference methods \cite{ShirokoffRosales2010} and meshfree finite difference method \cite{ZhouSeiboldShirokoffRosalesMeshfree2015}.
In \cite{ShirokoffRosales2010}, a second order finite difference scheme was proposed using a staggered spatial grid and a second order semi-implicit strategy in time where the viscous term was treated implicitly via Crank-Nicolson
and the pressure was treated explicitly using a second order Adams-Bashforth method.
Curved boundaries were embedded into a Cartesian grid, which could successfully handle irregular domains, but it was cumbersome to implement and did not generalize to arbitrary order.
Meshfree finite difference methods for the SR PPE reformulation were devised in \cite{ZhouSeiboldShirokoffRosalesMeshfree2015} as an alternative to allow for the  handling of irregular domains (without re-entrant corners). Specifically, \cite{ZhouSeiboldShirokoffRosalesMeshfree2015} devised
second order schemes with implicit-explicit time-stepping. Extensions to higher spatial order require larger stencils, thus lead to denser matrices and more costly computations.

In this paper we investigate a finite element discretization for a PPE reformulation with electric boundary conditions motivated by the SR PPE. The approach has the advantage of allowing for a systematic extension to higher spatial order.
We also use implicit-explicit (\imex{}) time-stepping, specifically \imex{} Runge-Kutta (RK) methods that decouple the velocity and the pressure solves.
We study the extent to which this yields desirable convergence and stability
properties --- i.e.\, better than second order in time, while avoiding a
parabolic time step restriction $\dt = O(\dx^2)$.
As a final note, we remark that the global existence of weak solutions to
the SR PPE reformulation, as well as to the Johnston and
Liu~\cite{JohnstonLiu2004} PPE reformulation, has been
proved~\cite{IgnatovaIyerKelliherPegoZarnescu2013} for no-slip boundary
conditions. The proofs provide appropriate function spaces for the weak
solutions, but they do not pursue finite element discretizations of the weak solutions.

This paper is organized as follows.
In \Srm\ref{sec:PPEforNSE}, we introduce a PPE reformulation motivated by the SR PPE of the Navier-Stokes equations. This reformulation uses non-standard boundary
conditions for the velocity --- i.e.\, electric boundary conditions (EBC). The
resulting vector Poisson-type problem with EBC, when solved with standard
nodal FEM, exhibits the Babu\u{s}ka paradox, which can be overcome by using
a mixed FEM formulation.
In \Srm\ref{sec:method_linear}, we present a numerical method for the \emph{linear} time-dependent problem (i.e.\, without the nonlinear advection term), based on mixed FEM and \imex{} RK time-stepping. We illustrate, via numerical tests, that the proposed method can achieve (at least) 3rd order in space and in time.
In \Srm\ref{sec:NonlinearAdvection}, we discuss how to extend the proposed
method to the nonlinear case. Numerical results are shown, for manufactured solutions, as well as for practical benchmark examples (lid-driven cavity, backward-facing step).

% --------------------------------------------------------------------------- 80
\section{Pressure Poisson Equation Reformulation of the
Navier-Stokes Equations} \label{sec:PPEforNSE}
% --------------------------------------------------------------------------- 80
In this section, we introduce a pressure Poisson equation reformulation of the Navier-Stokes equations with electric boundary conditions similar to the one proposed in~\cite{ShirokoffRosales2010}. Consider the time-dependent incompressible Navier-Stokes equations (NSE) in
a connected domain $\Omega\in\mathbb{R}^N$, where $N=2$ or $3$, with a
piece-wise smooth boundary $\partial\Omega$, for domains with Dirichlet boundary data,
\begin{subequations}\label{eq:NSE}
	%\left\{
\begin{alignat}{2}
\myvec{u}_t + \nl{u}{}{} & = \nu\Delta\myvec{u}-\nabla p + \myvec{f}
\quad && \text{in~}\Omega\times (0,T]\/, \label{eq:NSE1}\\
\DIV\myvec{u} & = 0 && \text{in~}\Omega\times (0,T]\/, \label{eq:NSE2}\\
\myvec{u}(\myvec{x},t) & = \myvec{g}(\myvec{x},t) &&
\text{on~}\partial\Omega\times [0,T]\/, \label{eq:NSE3} \\
\myvec{u}(\myvec{x},0) & = \myvec{u}_0(\myvec{x}) &&
\text{in~}\Omega\/, \label{eq:NSE4}
\end{alignat}%\right.
\end{subequations}
where $\nu > 0$ is the kinematic viscosity. Equation \eqref{eq:NSE1} follows
from the conservation of momentum, and \eqref{eq:NSE2} is conservation of
mass.
Furthermore, we impose the following compatibility conditions:
\begin{align}
&\textnormal{Continuity between the initial and the boundary conditions: } \myvec{u}_0(\myvec{x}) = \myvec{g}(\myvec{x},0) \text{ on~}\partial\Omega.\\
&\textnormal{Incompressibility of the initial condition: }\DIV\myvec{u}_0 = 0 \text{ in }\Omega.\\
&\textnormal{Zero net flux through the boundary: }\bdryint{\myvec{n}\cdot\myvec{g}} = 0.	\label{eq:NSE_zeroflux}
\end{align}

% --------------------------------------------------------------------------- 80
For the numerical solution of \eqref{eq:NSE} we will instead solve a PPE
reformulation, whose fundamental difference from previously proposed PPE
reformulations lies in the velocity boundary conditions:
incompressibility and the tangential flow are prescribed at the boundary.
Further, the normal velocity is enforced via a relaxation term in the
pressure equation. The PPE reformulation reads as follows
\begin{subequations}
\label{eq:SR-PPE_momentum}
\begin{alignat}{2}
  \myvec{u}_t + \nl{u}{}{} &= \nu\Delta\myvec{u}-\nabla P(\myvec{u}) +
  \myvec{f}\quad  &&\text{in } \Omega\times(0,T]\/,
  \label{eq:SR-PPE_momentum_1}\\
  \myvec{n} \times \myvec{u} &= \myvec{n}\times\myvec{g}
  &&\text{on }\partial\Omega\times[0,T]\/,
  \label{eq:SR-PPE_momentum_2}\\
  \nabla \cdot \myvec{u}  &= 0
  &&\text{on }\partial\Omega\times[0,T]\/. \label{eq:SR-PPE_momentum_3}
	\end{alignat}
\end{subequations}
Here $P(\myvec{u})$ is the solution to the pressure Poisson equation
(in a slight abuse of notation, we use both $p$ and $P(\myvec{u})$
for the pressure)
\begin{subequations}
\label{eq:SR-PPE_pressure}
\begin{alignat}{2}
  \Delta p &= \DIV(\myvec{f} - \nl{u}{}{})
  &&\text{in }\Omega\/,
  \label{eq:SR-PPE_pressure_1}\\
  \frac{\partial p}{\partial\myvec{n}}
   &= \myvec{n}\cdot \big(\myvec{f} -
      \myvec{g}_t-\nu\nabla\times\nabla\times\myvec{u}-\nl{u}{}{})\big) +
      \lambda \myvec{n}\cdot(\myvec{u}-\myvec{g}) \quad
   &&\text{on }\partial \Omega\/. \label{eq:SR-PPE_pressure_2}
\end{alignat}
\end{subequations}
The term $\lambda\myvec{n}\cdot(\myvec{u-g})$, where $\lambda > 0\/$
is a constant, is a relaxation term that guarantees that the normal
velocity condition at the boundary $\partial\Omega$ is exponentially
attracting --- see equation \eqref{eq:SR-normal}.

The tangential boundary conditions, together with the divergence-free boundary
condition in \eqref{eq:SR-PPE_momentum}, 
often appear in electrostatics as ``electric boundary conditions''
(EBC). We adopt this terminology in this paper. The PPE system (\ref{eq:SR-PPE_momentum}--\ref{eq:SR-PPE_pressure}) is almost identical to the one introduced in \cite{ShirokoffRosales2010}, the difference being here that $\myvec{n}\cdot\nabla\times\nabla\times \myvec{u}$ appears in the boundary condition for $p$ in lieu of
$\myvec{n}\cdot\Delta\myvec{u}$ (which appeared in \cite{ShirokoffRosales2010}). The choice of $\myvec{n}\cdot\nabla\times\nabla\times \myvec{u}$ in the pressure boundary condition is done to yield a simpler finite element discretization for the pressure than $\myvec{n} \cdot \Delta \myvec{u}$. Note that the pressure boundary condition $\myvec{n}\cdot\nabla\times\nabla\times \myvec{u}$ has appeared in several projection method and PPE formulations with Dirichlet boundary conditions for the velocity such as \cite{KarniadakisIsraeliOrszag1991, JohnstonLiu2004, LiuLiuPego2010}. We re-emphasize that our goal is to examine PPE schemes with EBC in the velocity.

In previous works, the PPE systems in \cite{ShirokoffRosales2010} and \cite{JohnstonLiu2004} were shown to be, for sufficiently smooth (up to the boundary) solutions $(\myvec{u},p)$, equivalent to the NSE \eqref{eq:NSE}. It is relatively straightforward to show that smooth solutions to the NSE solve (\ref{eq:SR-PPE_momentum}--\ref{eq:SR-PPE_pressure}) and hence the PPE system contains the NSE solutions.  We now show that solutions to the PPE system (\ref{eq:SR-PPE_momentum}--\ref{eq:SR-PPE_pressure}) solve the NSE --- the approach follows closely to \cite{ShirokoffRosales2010, JohnstonLiu2004} with a minor difference due to the combination of the EBC for the velocity and $\myvec{n}\cdot\nabla\times\nabla\times \myvec{u}$ boundary condition for the pressure. 

Assume that $(\myvec{u}, p)$ is a smooth solution to (\ref{eq:SR-PPE_momentum}--\ref{eq:SR-PPE_pressure}). Then the PPE reformulation recovers the incompressibility constraint: apply the divergence to the momentum equation \eqref{eq:SR-PPE_momentum_1} and
substitute into \eqref{eq:SR-PPE_pressure_1}. This yields the heat equation
for the divergence $\phi = \DIV\myvec{u}$, with homogeneous Dirichlet
boundary conditions due to \eqref{eq:SR-PPE_momentum_3}. That is:
\begin{equation}
	\label{eq:SR-heat}
	\begin{cases}
		\phi_t = \nu\Delta\phi & \text{in }\Omega\/,\\
		\phi = 0 & \text{on }\partial\Omega.
	\end{cases}
\end{equation}
Therefore, if $\phi(t=0) = \DIV\myvec{u}_0=0$, then $\phi=0$ for all time,
and $\myvec{u}$ is incompressible. If, due to numerical approximation
errors, the velocity field starts to depart from the $\DIV\myvec{u} = 0$
subspace, the heat equation dynamics ensure that $\myvec{u}$ is driven back
towards incompressibility. This property indicates that there is no need to
impose a discrete incompressibility principle in PPE reformulations, thus providing more flexibility in the design of numerical approximation methods.

Secondly, the PPE ensures that the normal velocity at the boundary, i.e.\, $\myvec{n}\cdot\myvec{u}= \myvec{n}\cdot\myvec{g}$, is enforced implicitly through the ordinary differential equation
\begin{equation}
\label{eq:SR-normal}
\begin{cases}
  \myvec{n}\cdot(\myvec{u}_t-\myvec{g}_t) = \nu \myvec{n}\cdot \nabla \phi -\lambda\myvec{n}\cdot(\myvec{u-g})
  & \text{on }\partial\Omega\times (0,T]\/,\\
  \myvec{n}\cdot(\myvec{u-g})|_{\partial\Omega} = 0 &\text{at } t=0,
\end{cases}
\end{equation}
which is obtained by evaluating the normal component of the momentum equation
\eqref{eq:SR-PPE_momentum_1} at the boundary, and using the pressure boundary
condition \eqref{eq:SR-PPE_pressure_2} along with the identity that $\Delta\myvec{u} = \nabla (\nabla \cdot \myvec{u}) - \nabla \times \nabla \times \myvec{u}$.  By virtue of equation \eqref{eq:SR-heat}, we have that $\phi = 0$, so that \eqref{eq:SR-normal} reduces to a simple ODE at every point $\myvec{x} \in \partial \Omega$: $\dot{\alpha} = -\lambda \alpha$ where $\alpha = \myvec{n}\cdot(\myvec{u}-\myvec{g})$. Hence, 
$\myvec{n}\cdot\myvec{u}=\myvec{n}\cdot\myvec{g}$ for all time. In the
presence of numerical approximation errors, the term
$-\lambda\myvec{n}\cdot(\myvec{u}-\myvec{g})$ adds an exponential decay to
the error in the normal velocity on the domain boundary.

Note that in the absence of the $\lambda$-term, the exact solution to
equation \eqref{eq:SR-normal} still satisfies
$\myvec{n}\cdot(\myvec{u}-\myvec{g})=0$. However, in practice, numerical
errors may result in a drift of the normal velocity \cite{ShirokoffRosales2010},
which requires stabilization via the $\lambda$-term.
Moreover, for the steady state system, the condition $\lambda>0$ is required
for a unique solution.

% --------------------------------------------------------------------------- 80
\begin{remark}[Solvability condition for the pressure Poisson equation]
\label{rem:PPE_solvability}
A solvability condition is required in order for the pressure Poisson
equation \eqref{eq:SR-PPE_pressure} to have a solution, that is:
\begin{equation}\label{eq:solvability}
\int_{\Omega} -\nu\nabla \cdot (\nabla \times \nabla \times \myvec{u})+\lambda\DIV\myvec{u}\intd{V} -
\bdryint{\myvec{n}\cdot\myvec{g}_t + \lambda\myvec{n}\cdot\myvec{g}}=0\/.
\end{equation}
The incompressibility condition $\DIV\myvec{u}=0$ and the zero net flux
condition \eqref{eq:NSE_zeroflux}
ensure that both the volume integral and
the boundary integral in equation \eqref{eq:solvability} vanish. However,
numerical approximation errors in the discrete Poisson equation may result in schemes
that do not exactly satisfy the discrete version of the solvability
condition \eqref{eq:solvability}.
	Whenever this occurs, the discrete Poisson equation
	must be solved in the least squares sense.
This is achieved by formulating an augmented system that projects the right
hand side of the pressure Poisson problem onto one for which the solvability
condition is satisfied, see \Srm\ref{Sec:AppendixPoisson}.
\end{remark}
% --------------------------------------------------------------------------- 80
%
% --------------------------------------------------------------------------- 80
\begin{remark}\label{PPE:richness}
There is a rich variety of possible PPE reformulations.
At the continuum level, one can add $\DIV\myvec{u}$ anywhere (because $\DIV\myvec{u}=0$), as long as the resulting systems are well-posed and are equivalent to the original problem. However, in the presence of approximation errors, one generally has $\DIV\myvec{u}_h\neq 0$, thus adding $\DIV\myvec{u}$ will lead to different numerical schemes.
For example, Henshaw and Petersson \cite{HenshawPetersson2003} add a divergence damping term $\delta\,\DIV\myvec{u}$, with $\delta\geq 0$, to the pressure Poisson equation to obtain
\begin{equation*}
 \Delta p = \DIV(\myvec{f}-\nl{\myvec{u}}{}{}) + \delta \DIV\myvec{u}.
\end{equation*}
With this, the divergence satisfies the PDE
$\phi_t = \Delta\phi -\delta\phi$. In the discretized case, where the divergence
is not exactly zero, the damping term adds an exponential decay that can further help keeping the discrete divergence small.
\end{remark}
% --------------------------------------------------------------------------- 80

% --------------------------------------------------------------------------- 80
\section{Numerical Method for the Time-dependent Stokes Problem}
\label{sec:method_linear}
% --------------------------------------------------------------------------- 80
%
In this section, we present a numerical method for the PPE reformulation
(\ref{eq:SR-PPE_momentum}--\ref{eq:SR-PPE_pressure}) of the time-dependent Stokes equation
\begin{equation}
\label{eq:linear_SR-PPE_u}
  \myvec{u}_t = \nu\Delta\myvec{u}-\nabla P(\myvec{u}) +
  \myvec{f}\quad \text{ in } \Omega\times(0,T]\/,\quad 
  \myvec{n} \times \myvec{u} = \myvec{n}\times\myvec{g}
  \mbox{ and }  \nabla \cdot \myvec{u}  = 0
  \quad\text{on }\partial\Omega\times[0,T]\/,
\end{equation}
where $P(\myvec{u})$ solves the pressure Poisson equation
\begin{equation}
\label{eq:linear_SR-PPE_p}
  \Delta p = \DIV\myvec{f}
  \quad\text{ in }\Omega\/,\quad
  \frac{\partial p}{\partial\myvec{n}}
   = \myvec{n}\cdot (\myvec{f} -
      \myvec{g}_t-\nu\nabla\times\nabla\times\myvec{u}) +
      \lambda \myvec{n}\cdot(\myvec{u}-\myvec{g}) \quad
   \text{ on }\partial \Omega\/. 
\end{equation}
To handle irregular domains we
adopt a mixed finite element method (FEM) for the spatial discretizations
of $\myvec{u}$ (see \Srm\ref{ssec:method_linear_velocity}), and nodal FEM for
$p$ (see \Srm\ref{ssec:method_linear_pressure}).
Note that the motivation for adopting a mixed FEM for $\myvec{u}$ is to address structural issues (discussed below) that arise from the electric boundary conditions (EBC) in the momentum equation.
Unlike mixed approaches for the Stokes/Navier-Stokes equations in which the elements approximating $\myvec{u}$ and $p$ need to satisfy a discrete inf-sup condition \cite{FortinBrezzi1991,Layton2008},
the mixed formulation employed here is solely for the velocity, and the nodal elements for the pressure do not need to satisfy an inf-sup condition with the velocity elements.
For the time evolution (see \Srm\ref{ssec:method_linear_time}), we adopt an
implicit-explicit (\imex{}) Runge-Kutta (RK)
scheme that (due to the PPE formulation) leads to natural approaches
for decoupling the velocity $\myvec{u}$ from the pressure $p$. That is,
we treat the viscous term $\Delta \myvec{u}$ implicitly and the pressure
term $\nabla p$ explicitly.
We carry out convergence studies via the method of manufactured
solutions in \Srm\ref{ssec:linear_NumerialResults}.

% --------------------------------------------------------------------------- 80
\subsection{Discretization of the Velocity via Mixed Finite Elements}
\label{ssec:method_linear_velocity}
% --------------------------------------------------------------------------- 80
%
In this section we outline the spatial discretization of the velocity via
mixed finite elements \cite{FortinBrezzi1991}.  The choice of a mixed FEM
over other element choices (such as nodal FEM) is due to the
EBC in the momentum equation. To be precise, nodal FEM for problems
involving electric boundary conditions may converge to the wrong
solution, see Remark~\ref{rmk:failurecvg}. In contrast, mixed FE provide a natural way to handle the EBC.
% --------------------------------------------------------------------------- 80
%
% --------------------------------------------------------------------------- 80

The mixed formulation presented in this paper is applied to the discretization of the momentum equation \eqref{eq:SR-PPE_momentum} only, and introduced to handle the EBC by introducing the vorticity $\sigma=\CURL\myvec{u}$ as a new variable. Hence, there is no inf-sup condition for $(\myvec{u}, p)$, but instead, an inf-sup condition for the velocity and vorticity $(\myvec{u}, \sigma)$.  Using the vector identity
$\Delta\myvec{u} = \nabla(\DIV\myvec{u})-\CURL\CURL\myvec{u}$, and introducing
the new variable
$\sigma = \CURL\myvec{u}$ (cf.~\cite{ArnoldFalkGopalakrishnan2012,ArnoldChen2017}),
equation \eqref{eq:linear_SR-PPE_u} in the PPE
reformulation can be recast as:
% --------------------------------------------------------------------------- 80
\begin{subequations}
	\begin{alignat}{2}
		\sigma &= \CURL\myvec{u} &&\text{in }\Omega\label{eq:mixed_SR-PPE1},\\
		\myvec{u}_t &= \nu(\nabla(\DIV\myvec{u})-\CURL\sigma)-\nabla P(\myvec{u})+\myvec{f}\quad &&\text{in }\Omega,\label{eq:mixed_SR-PPE2}\\
		\myvec{n}\times\myvec{u} &= \myvec{n}\times\myvec{g}&&\text{on }\partial\Omega,\label{eq:mixed_SR-PPE3}\\
		\DIV\myvec{u} &= 0&&\text{on }\partial\Omega\label{eq:mixed_SR-PPE4}.
	\end{alignat}
	\label{eq:mixed_SR-PPE}
\end{subequations}
% --------------------------------------------------------------------------- 80
To obtain the weak formulation of the equations \eqref{eq:mixed_SR-PPE}, we
use the spaces
\begin{align*}
H(\tcurl;\Omega) &= \left\{\myvec{u}\in L^2(\Omega)^N:\,\CURL\myvec{u}\in L^2(\Omega)^{N'}\right\},\\
H(\tdiv;\Omega)& = \left\{\myvec{u}\in L^2(\Omega)^N:\,\DIV\myvec{u}\in L^2(\Omega)\right\},
\end{align*}
with $N'=1$ for $N=2$ and $N'=3$ for $N=3$.
We then multiply \eqref{eq:mixed_SR-PPE1} by a test function $\tau$ and
\eqref{eq:mixed_SR-PPE2} by a test function $\myvec{v}$, and then apply the
integral identities
% --------------------------------------------------------------------------- 80
\begin{align*}
	\inprod{\CURL\myvec{u}}{\tau} &=  \inprod{\myvec{u}}{\CURL\tau}+\bdryint{\tau\cdot(\myvec{n}\times\myvec{u})},\\
	\inprod{\CURL\sigma-\nabla(\DIV\myvec{u})}{\myvec{v}}& = \inprod{\CURL\sigma}{\myvec{v}} + \inprod{\DIV\myvec{u}}{\DIV\myvec{v}} - \bdryint{(\DIV\myvec{u})(\myvec{v}\cdot\myvec{n})},
\end{align*}
% --------------------------------------------------------------------------- 80
along with the tangential boundary condition \eqref{eq:mixed_SR-PPE3} and the divergence-free boundary condition \eqref{eq:mixed_SR-PPE4}. This procedure yields the following mixed formulation for \eqref{eq:mixed_SR-PPE}:
Find $\sigma\in X$, $\myvec{u}\in H(\tdiv;\Omega)$ such that
% --------------------------------------------------------------------------- 80
\begin{subequations}
	\label{eq:mixed_NSEPPE}
	\begin{alignat}{2}
		\inprod{\sigma}{\tau} &- \inprod{\myvec{u}}{\CURL\tau} = \bdryint{\tau(\myvec{n}\times\myvec{g})} &&\forall\tau\in X,\\
		\inprod{\myvec{u}_t}{\myvec{v}}&+\nu \inprod{\CURL\sigma}{\myvec{v}} + \nu\inprod{\DIV\myvec{u}}{\DIV\myvec{v}} = \inprod{\myvec{f}-\nabla P(\myvec{u})}{\myvec{v}}\quad &&\forall\myvec{v}\in H(\tdiv;\Omega).
	\end{alignat}
\end{subequations}
% --------------------------------------------------------------------------- 80
Here $X = H^1(\Omega)$ when $N=2$ and $X = H(\tcurl;\Omega)$ when $N=3$.
Notice that both the tangential boundary condition and the divergence
boundary condition from \eqref{eq:SR-PPE_momentum}
appear in the weak form \eqref{eq:mixed_NSEPPE} as natural boundary
conditions (the spaces of functions $X$ and $H(\tdiv; \Omega)$ do not
enforce the boundary conditions).
% --------------------------------------------------------------------------- 80

% --------------------------------------------------------------------------- 80
Let $\Sigma^h$ and $V^h$ be finite dimensional subspaces for $X$ and
$H(\tdiv;\Omega)$ respectively. The semi-discrete mixed formulation of
\eqref{eq:mixed_NSEPPE} is as follows:
Find $(\sigma_h,\myvec{u}_h)\in \Sigma^h\times V^h$ such that
$\myvec{u}_h(t=0)={\myvec{u}_0}_h$ and
\begin{subequations}
	\label{eq:mixed_SR-NSE_semi}
	\begin{alignat}{2}
		\inprod{\sigma_h}{\tau_h} &- \inprod{\myvec{u}_h}{\CURL\tau_h} = \bdryint{\tau_h(\myvec{n}\times\myvec{g})}\qquad && \quad\forall \tau\in\Sigma^h\\
		\inprod{(\myvec{u}_h)_t}{\myvec{v}_h} &= -\nu \inprod{\DIV\myvec{u}_h}{\DIV\myvec{v}_h}-\nu\inprod{\CURL\sigma_h}{\myvec{v}_h}
		+ \inprod{\myvec{f}-\nabla P(\myvec{u}_h)}{\myvec{v}_h} && \quad\forall \myvec{v}_h\in V^h
	\end{alignat}
\end{subequations}
where ${\myvec{u}_0}_h$ is the projection of the initial condition
$\myvec{u}_0$ onto the space $V^h$, and $P(\myvec{u}_h)$ is the solution to
the discretized pressure Poisson problem using the velocity approximation
(see \Srm\ref{ssec:method_linear_pressure}).

There are different stable pairs of finite elements for the spaces $\Sigma^h$ and $V^h$
\cite{ArnoldFalkGopalakrishnan2012,ArnoldFalkWinther2006,ArnoldFalkWinther2010}.
Here in the case of $N=2$, we choose nodal finite elements
(Lagrange finite elements $P_r$) of degree $r\geq 1$ for $\sigma_h$ and
Raviart-Thomas elements ($RT_{r-1}$) of the same degree for the vector field $\myvec{u}_h$.
That is
\begin{equation*}
\Sigma^h\times V^h = P_r\times RT_{r-1}\quad\mbox{for~}r\geq 1.
\end{equation*}
In the 3-dimensional case, the corresponding space for $\sigma$ is $H(\tcurl;\Omega)$, and N\'{e}d\'{e}lec elements ($NED_{r}^1$) are used.

% --------------------------------------------------------------------------- 80
% 
% --------------------------------------------------------------------------- 80

	\begin{remark}\label{rmk:failurecvg}
		(Failure of nodal FEM)
		Nodal FEM may fail to converge to the true solution for EBC problems, for two different reasons:
		(i)~the Babu\u{s}ka paradox \cite{Babuska1961,BabuskaPitkaranta1990}, and
		(ii)~the inability to approximate singularities in the solution, such as
		those caused by re-entrant corners in the domain.
		The Babu\u{s}ka paradox \cite{Babuska1961,BabuskaPitkaranta1990} occurs when FEM solutions on polygonal approximations of domains with curved boundaries converge to a function (as the mesh size goes to zero) that is not the solution of the underlying continuum problem.
		For instance, the vector Laplace problem (with $\DIV\myvec{f} = 0$)
		\begin{equation}\label{eq:vpe}
		\Delta\myvec{u} = \myvec{f}\quad\mbox{in~}\Omega\/,
		\quad\mbox{with b.c.~}\quad
		\DIV\myvec{u}=0\quad\mbox{and}\quad
		\myvec{n}\times\myvec{u}=0\quad\mbox{on~}\partial\Omega\/,
		\end{equation}
		has the weak formulation: Given $\myvec{f}\in L^2(\Omega)^N$, $\DIV\myvec{f}=0$, find
		$\myvec{u}\in V$, such that for every $\myvec{v}\in V$
		\begin{align}\label{eq:VP1}
		\inprod{\CURL\myvec{u}}{\CURL\myvec{v}}+
		\inprod{\DIV\myvec{u}}{\DIV\myvec{v}} = \inprod{\myvec{f}}{\myvec{v}}.
		\end{align}
		Nodal FEM approximations based on $V = H^1_{0t}(\Omega)^N = \{\myvec{u}\in H^1(\Omega)^N: \myvec{n}\times\myvec{u}|_{\partial\Omega}=0\}$ exhibit the Babu\u{s}ka paradox  and may not converge to the solution of \eqref{eq:vpe} \cite{ZhouSeiboldShirokoffRosalesMeshfree2015}.
		
		The second issue with nodal elements is an approximation theory result.
		The weak formulation \eqref{eq:VP1} is uniquely solvable for the following two choices of space: $V=H^1_{0t}(\Omega)^N$ or $V = H_{0}(\tcurl;\Omega)\cap H(\tdiv;\Omega) = \{\myvec{u}\in L^2(\Omega)^N: \DIV\myvec{u}\in L^2(\Omega)^N, \CURL\myvec{u}\in L^2(\Omega)^{N'}, \myvec{n}\times\myvec{u}|_{\partial\Omega}=0 \}$. 
		When the domain $\Omega$ is convex or has a globally $C^2$ boundary, then $H^1_{0t}(\Omega)^N = H_{0}(\tcurl;\Omega)\cap H(\tdiv;\Omega)$ and the solution is divergence-free.
		When $\Omega$ has re-entrant corners, we have $H^1_{0t}(\Omega)^N \subsetneq H_{0}(\tcurl;\Omega)\cap H(\tdiv;\Omega)$ and the divergence-free solution may not in $H^1_{0t}(\Omega)^N$\cite{KangroNicholaides1999}. Hence, in general $H^1$ nodal FEM is not guaranteed to converge to the divergence-free solution.
	\end{remark}

% --------------------------------------------------------------------------- 80
\subsection{Discretization of the Pressure}
\label{ssec:method_linear_pressure}
% ======================================================================= @#@ 80
The pressure $p$ satisfies a Poisson problem with Neumann boundary condition
\eqref{eq:linear_SR-PPE_p}. Therefore, given the velocity, standard
nodal-based finite elements can be used to discretize the pressure Poisson equation. Following the usual procedure for deriving weak formulations, we multiply equation \eqref{eq:linear_SR-PPE_p} by a test function $q\in H^1(\Omega)$ and apply integration by parts to obtain
\begin{equation*}
	\inprod{\nabla p}{\nabla q} =  \inprod{\myvec{f}}{\nabla q} -
	\nu\bdryint{(\myvec{n}\cdot\nabla \times \nabla \times\myvec{u})q}
	+\lambda\bdryint{\myvec{n}\cdot(\myvec{u}-\myvec{g})q}
	-\bdryint{(\myvec{n}\cdot \myvec{g}_t)q}.
\end{equation*}
We further substitute $\sigma = \nabla \times \myvec{u}$, which gives rise to two possible different weak formulations for the pressure Poisson equation:
\begin{enumerate}
	\item Choosing the boundary integral that involves $\myvec{n}\cdot(\CURL\sigma)$ yields the first weak formulation
	\begin{equation}
		\inprod{\nabla p}{\nabla q} =  \inprod{\myvec{f}}{\nabla q}
		- \nu\bdryint{\myvec{n}\cdot(\CURL\sigma)q}
		+\lambda\bdryint{\myvec{n}\cdot(\myvec{u}-\myvec{g})q}
		-\bdryint{(\myvec{n}\cdot \myvec{g}_t)q},\quad\forall q\in H^1(\Omega).
		\label{eq:pressure_weak1}
	\end{equation}
	\item Choosing the volume integral of $(\CURL\sigma)\cdot\nabla q$ gives the second weak formulation
	\begin{equation}
		\inprod{\nabla p}{\nabla q} = \inprod{\myvec{f}}{\nabla q}
		- \nu\inprod{\CURL\sigma}{\nabla q}
		+\lambda\bdryint{\myvec{n}\cdot(\myvec{u}-\myvec{g})q}
		-\bdryint{(\myvec{n}\cdot \myvec{g}_t)q},\quad\forall q\in H^1(\Omega).
		\label{eq:pressure_weak2}
	\end{equation}
\end{enumerate}
In numerical tests both choices return almost identical results. Therefore in
this paper we show only the results of numerical experiments conducted with
the first formulation, \eqref{eq:pressure_weak1}.

% --------------------------------------------------------------------------- 80
\subsection{Time-stepping via \imex{} Schemes}
\label{ssec:method_linear_time}
% --------------------------------------------------------------------------- 80
One advantage of the PPE reformulation is that the pressure appears as a global function of the velocity in the momentum equation. This enables conceptually straightforward implicit-explicit (\imex{}) time discretizations of the PPE system that treat the pressure explicitly and viscosity implicitly. An implicit treatment of the viscosity term is desirable for low to moderate Reynolds numbers as to avoid a parabolic time step stability restriction. In addition, an explicit treatment of the pressure is desirable to avoid solving large coupled systems involving $(p, \myvec{u})$.

\imex{} schemes are based on an additive splitting of an ODE that take the following form:
\begin{equation}\label{eq:imex_ode}
 \dfrac{\text{d} u}{\text{d}t} = f(u)+g(u)\/.	
\end{equation}
Popular \imex{} methods are linear multistep \imex{} schemes \cite{AscherRuuthWetton1995} and \imex{} Runge-Kutta (RK) methods \cite{AscherRuuthSpiteri1997}. Here we focus on \imex{} RK schemes because they have less restrictive stability properties than \imex{} multistep methods.
We consider \imex{} RK schemes that combine two different Runge-Kutta schemes: an explicit RK (ERK) method for $f(u)$, and a diagonally implicit RK (DIRK) method for $g(u)$.
Let $A\in\mathbb{R}^{s\times s}$, $\myvec{b},\myvec{c}\in\mathbb{R}^s$ be the
coefficients of an $s$-stage DIRK scheme, and
$\hat{A}\in\mathbb{R}^{(s+1)\times (s+1)}$,
$\hat{\myvec{b}},\,\hat{\myvec{c}}\in\mathbb{R}^{s+1}$
be the coefficients of an $(s+1)$-stage ERK scheme with
$\hat{\myvec{c}}^{T} = (0,\myvec{c}^{T})$ in the Butcher
notation \cite{Butcher1963}, i.e.\,
\begin{align*}
	\begin{array}{c|c}
		\myvec{c} & A \\
		\hline
		& \myvec{b}^T
	\end{array}
	=
	\begin{array}{c|c c c c }
		c_1 & a_{11} & & & \\
		c_2 & a_{21} & a_{22} & & \\
		\vdots & \vdots & \vdots & \ddots &   \\
		c_s & a_{s1} & a_{s2} & \cdots& a_{ss} \\
		\hline
		& b_1 & b_2 & \cdots & b_s
	\end{array}
	\quad\quad
	\begin{array}{c|c}
		\hat{\myvec{c}} & \hat{A}\\
		\hline
		& \hat{\myvec{b}}^T
	\end{array}
	=
	\begin{array}{c|c c c c c}
		0 & 0 & & & & \\
		c_1  & \hat{a}_{21}& 0 & & & \\
		c_2 & \hat{a}_{31} & \hat{a}_{32}& 0 & & \\
		\vdots & \vdots & \vdots & \vdots & \ddots &   \\
		c_s & \hat{a}_{s+1,1} & \hat{a}_{s+1,2} & \hat{a}_{s+1,3} & \cdots& 0 \\
		\hline
		&\hat{b}_1 & \hat{b}_2 & \hat{b}_3 & \cdots & \hat{b}_{s+1}
	\end{array}
\end{align*}
One step of an \imex{} scheme from time $t_n$ to $t_{n+1}=t_n+\dt$ for the
splitting \eqref{eq:imex_ode} can be written as follows:
\begin{align*}
	u^{(i)} & = u^{n} +\dt\sum_{j=1}^i a_{ij}g(u^{(j)}) +\dt\sum_{j=1}^{i} \hat{a}_{i+1,j}f(u^{(j-1)})\quad\mbox{for}\quad i = 1,\dots,s \\
	u^{n+1} & = u^n + \dt\sum_{i=1}^s b_j g(u^{(j)}) + \dt\sum_{j=1}^{s+1}\hat{b}_j f(u^{(j-1)})
\end{align*}
where $u^{(0)} = u^n$.

% --------------------------------------------------------------------------- 80
% --------------------------------------------------------------------------- 80
\vspace{1em}
\noindent	
\textbf{Stability of \imex{} RK Schemes for the PPE Reformulation.} IMEX schemes applied to the PPE system (\ref{eq:SR-PPE_momentum}--\ref{eq:SR-PPE_pressure}) may encounter subtle stability issues, and addressing these issues will guide our choice of IMEX time stepping coefficients.  First we will examine the numerical stability of an \imex{} RK scheme applied to a scalar ODE model for the PPE system (\ref{eq:linear_SR-PPE_u}--\ref{eq:linear_SR-PPE_p}). The model problem will then yield a stability criterion, which will dictate our choice of \imex{} coefficients.

Using the vector identity $\Delta = \nabla(\nabla \cdot) - \nabla \times \nabla \times$, and introducing $\mathcal{A} := -\nabla \times \nabla \times$, the PPE system (\ref{eq:linear_SR-PPE_u}--\ref{eq:linear_SR-PPE_p}) can be (exactly) recast (with $\nu=1$, when $\myvec{f} = 0$, $\myvec{g} = 0$) as the following non-local evolution
\begin{equation}
	\label{eq:SRPPEprojection}
	\myvec{u}_t = \nabla(\nabla\cdot\myvec{u}) + \mathcal{A}\myvec{u} 
	-\mathcal{P}(\mathcal{A}\myvec{u}) - \lambda \mathcal{P}\myvec{u}
	\quad\mbox{in }\Omega\quad
	\mbox{with b.c. }\myvec{n}\times\myvec{u}=0,\;\;\DIV\myvec{u}=0\quad\mbox{on }\partial\Omega.
\end{equation}

Here the operator $\mathcal{P}$ is defined by
\begin{equation}
	\mathcal{P}\myvec{w} = \nabla p\/,
\end{equation}
where $p$ solves the Poisson equation
\begin{equation}
\label{eq:SRPPEprojection_poisson}
\Delta p = 0\quad\mbox{in }\Omega\/,
\quad\mbox{with b.c. }
\frac{\partial p}{\partial\myvec{n}} = \myvec{n}\cdot\myvec{w} - \bar{w}
\quad\mbox{on~}\partial\Omega,\quad
\mbox{where~}
\bar{w}:=\frac{1}{|\partial\Omega|}\bdryint{\myvec{n}\cdot\myvec{w}}.
\end{equation}
By definition, $\mathcal{P}\myvec{w}$ generates a divergence-free field, with
normal component $\myvec{n}\cdot\myvec{w}-\bar{w}$ at the boundary. It
follows that $\mathcal{P}^2 = \mathcal{P}$, i.e.\, $\mathcal{P}$ is a projection.

The IMEX discretization of \eqref{eq:SRPPEprojection} then treats $\nabla(\nabla\cdot\myvec{u}) + \mathcal{A}\myvec{u}$ implicitly and $-\mathcal{P}(\mathcal{A}\myvec{u}) - \lambda \mathcal{P}\myvec{u}$ explicitly. Due to the structure of the projection, $\mathcal{P}(\mathcal{A}\myvec{u})\/$ is
stiff since $\mathcal{A}$ has two spatial derivatives 
(i.e.\, may incur a parabolic time step restriction). 
Thus both the implicit term (which includes $\mathcal{A}\myvec{u}$) and the explicit term ($\nabla p$, which includes $\mathcal{P}\mathcal{A}\myvec{u}$) are stiff. To investigate the stability of \imex{} RK schemes for
solving the PPE reformulation, we consider a scalar model problem for equation \eqref{eq:SRPPEprojection}
\begin{equation}\label{eq:testproblem}
	u_t = -\gamma u+\mu u\/,
\end{equation}
where $-\gamma u$ is treated implicitly and $\mu u$ explicitly. The model equation \eqref{eq:testproblem} is frequently used to understand \imex{} stability. However, commonly one considers $\gamma>0$ real and $\mu$ purely imaginary, e.g.\, to understand time-stepping for advection-diffusion problems \cite{AscherRuuthWetton1995,AscherRuuthSpiteri1997}.
In contrast, we take \emph{both} $\gamma$ and $\mu$ to be real. Specifically, we focus on the interplay of $\mathcal{A}\myvec{u}$ and $\mathcal{P}\mathcal{A}\myvec{u}$ and ignore $\nabla(\nabla\cdot\myvec{u})$ since it is stabilizing, and $\lambda \mathcal{P}\myvec{u}$ since it is nonstiff. We then model $\gamma$ as an eigenvalue of $-\mathcal{A}\myvec{u}$ and $\mu$ as an eigenvalue of $\mathcal{P}\mathcal{A}\myvec{u}$. Since $-\mathcal{A}$ is a symmetric positive semi-definite operator (in the $L^2$ inner product on vector fields with EBC), and $\mathcal{P}$ is a projection, we restrict the parameter values by $0<\mu\leq\gamma$\/. The situation $(\gamma,\mu)\in\mathbb{R}^2$ was studied in \cite{Koto2008}, to understand the stability properties for certain 1st and 2nd order \imex{} RK schemes in the context of reaction-diffusion equations. However, both the stiffness behavior and the important wedge property introduced below (see Definition~\ref{def:wedge_property}) were not discussed in \cite{Koto2008}.

Let us pad the $s$-stage implicit Runge-Kutta scheme with zeros to obtain
the tableau
\begin{equation}
	\begin{array}{c|c}
		\tilde{\myvec{c}} & \tilde{A} \\
		\hline
		& \tilde{\myvec{b}}^T
	\end{array}
	=
	\begin{array}{c|c c c }
		0     		& 0 & 0  \\
		\myvec{c} & 0 & A  \\
		\hline
		& 0 & \myvec{b}^T
	\end{array} \/,
\end{equation}
and introduce the vector $\myvec{U}^{n+1} = (u^n, u^{n+1}_1,\dots,u^{n+1}_s)^T$ that includes the $s$ intermediate stage solutions and the approximation from previous time step. Then one step of the \imex{} RK scheme applied to the model problem \eqref{eq:testproblem} can be expressed as
\begin{align}
	\label{eq:testproblem_imex_1}
	\myvec{U}^{n+1} &= u^n\myvec{e} - \gamma\dt\tilde{A}\myvec{U}^{n+1} + \mu\dt\hat{A}\myvec{U}^{n+1}\/,\\
	\label{eq:testproblem_imex_2}
	u^{n+1} &= u^n - \gamma\dt\tilde{\myvec{b}}^T\myvec{U}^{n+1} + \mu\dt\hat{\myvec{b}}^T\myvec{U}^{n+1}\/,
\end{align}
where $\myvec{e}$ is an $(s+1)$-vector of all ones. The above expressions yield
\begin{equation}
	u^{n+1} = R(\gamma\dt,\mu\dt) u^n\/,
\end{equation}
where
\begin{align*}
	R(\alpha,\beta) &= 1+(-\alpha\tilde{\myvec{b}}^T + \beta\hat{\myvec{b}}^T)(I+\alpha\tilde{A} -\beta\hat{A})^{-1}\myvec{e}\\
	& = \frac{\det(I+\alpha\tilde{A}-\beta\hat{A}-\alpha\myvec{e}\tilde{\myvec{b}}^T +\beta\myvec{e}\hat{\myvec{b}}^T)}{\det(I+\alpha\tilde{A}-\beta\hat{A})}.
\end{align*}
is the \imex{} RK analogue of the stability function for Runge-Kutta schemes
applied to the test problem $y'=\lambda y$.
It is now natural to define the \emph{stability region} as
\begin{equation}
	S = \{(\alpha,\beta):\, |R(\alpha,\beta)| \leq 1\}.
\end{equation}
For a given \imex{} RK scheme, we can plot the stability region in the $(\alpha,\beta)$-plane, where $\alpha = \gamma\dt$ and $\beta = \mu\dt$. In general, when both $\gamma$ and $\mu$ are in $\mathbb{C}$, the stability region $S$ is a subset in $\mathbb{C}^2$.
For the model problem \eqref{eq:testproblem}, it suffices to look at the cross section of the region $S$ in the first quadrant of $\mathbb{R}^2$, which can be easily visualized.
\begin{definition}\label{def:wedge_property}
An \imex{} RK scheme is said to possess the wedge property, if the stability region $S$ contains the wedge $\{(\alpha,\beta):\, 0<\beta\leq\alpha\}$.
In other words, the scheme is unconditionally stable for the model problem \eqref{eq:testproblem} with any $0<\mu\le\gamma$.
\end{definition}

\begin{figure}[htb!]
 \includegraphics[width=0.32\textwidth]{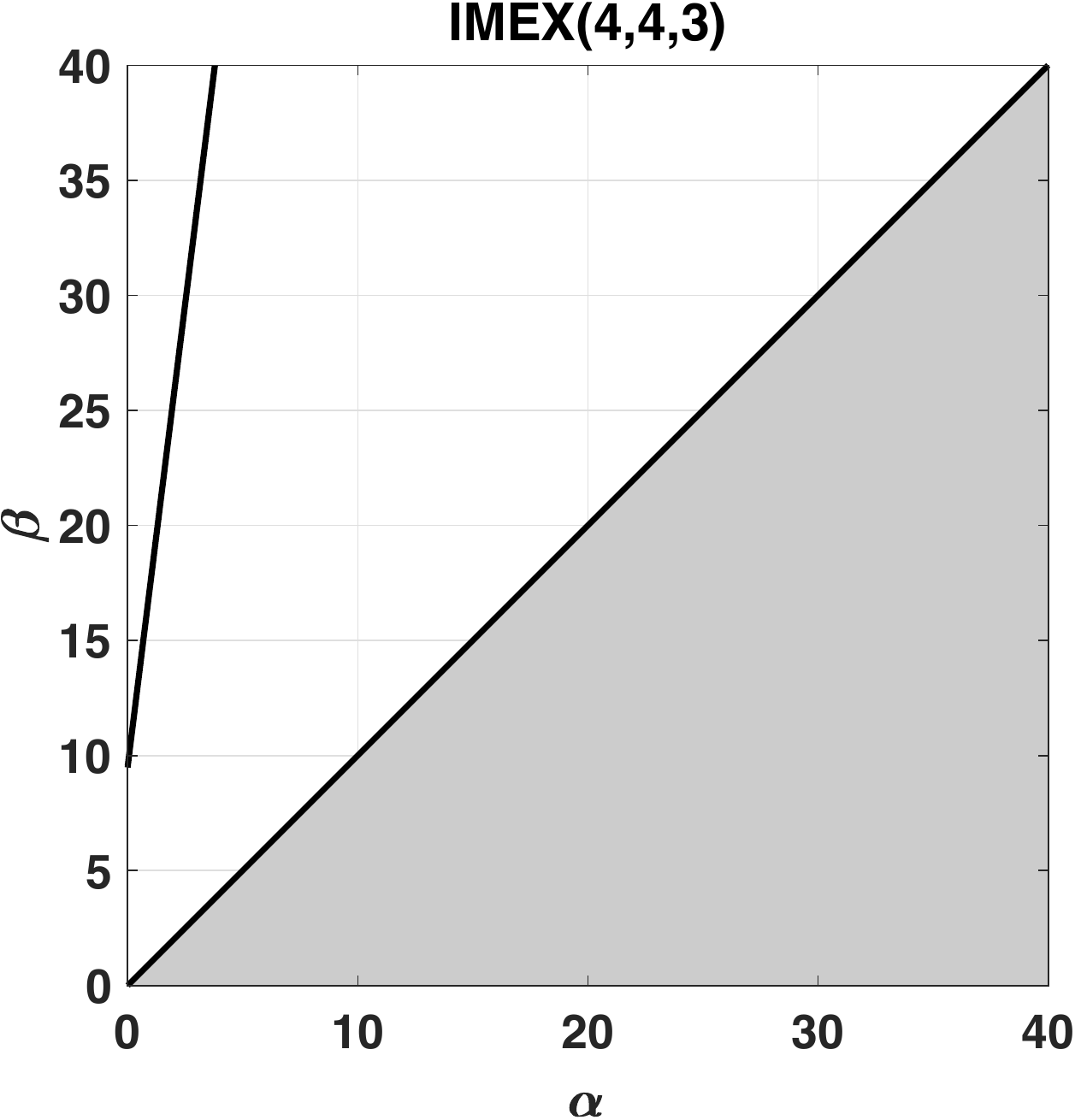}
 \includegraphics[width=0.32\textwidth]{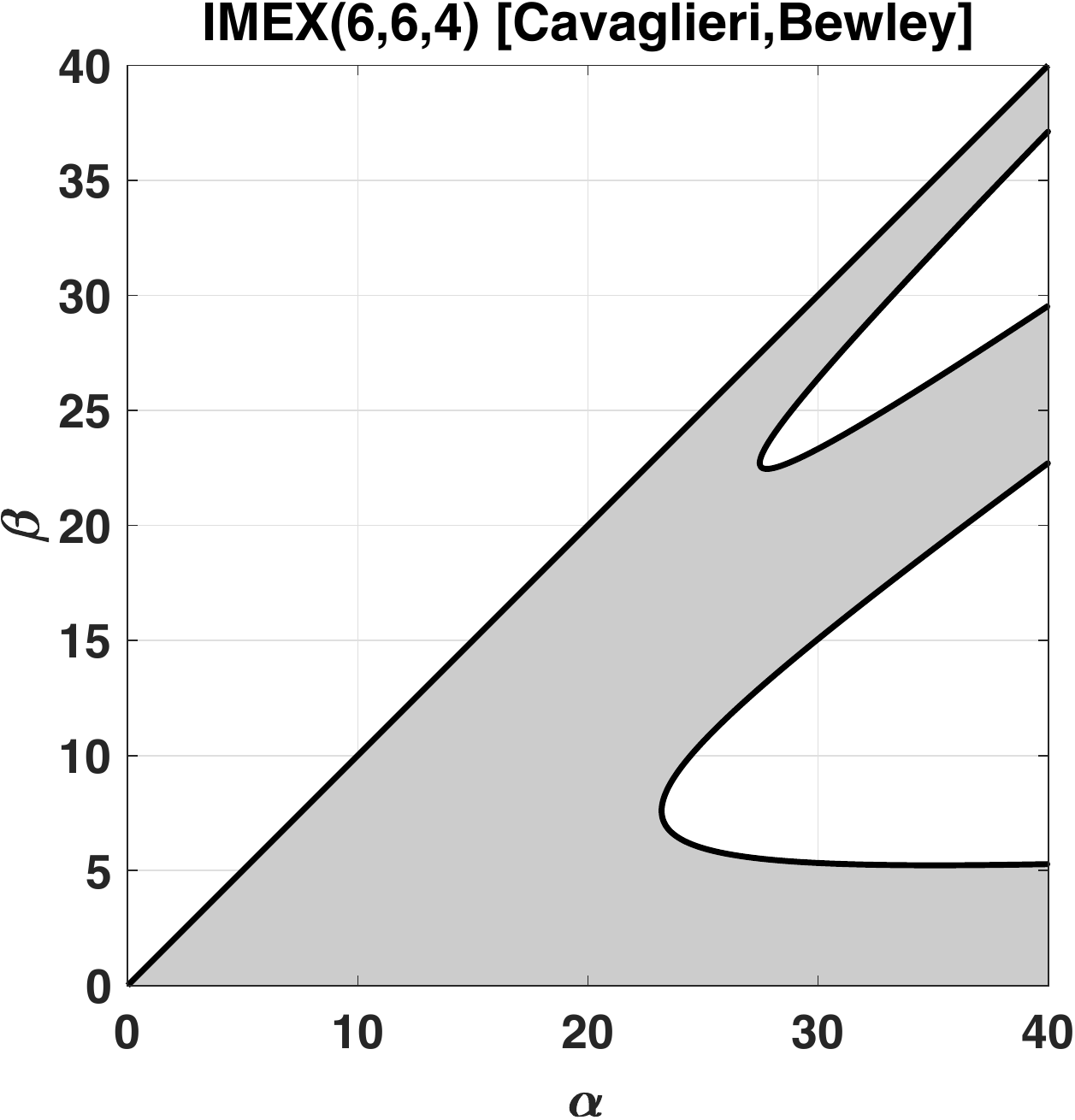}
 \includegraphics[width=0.32\textwidth]{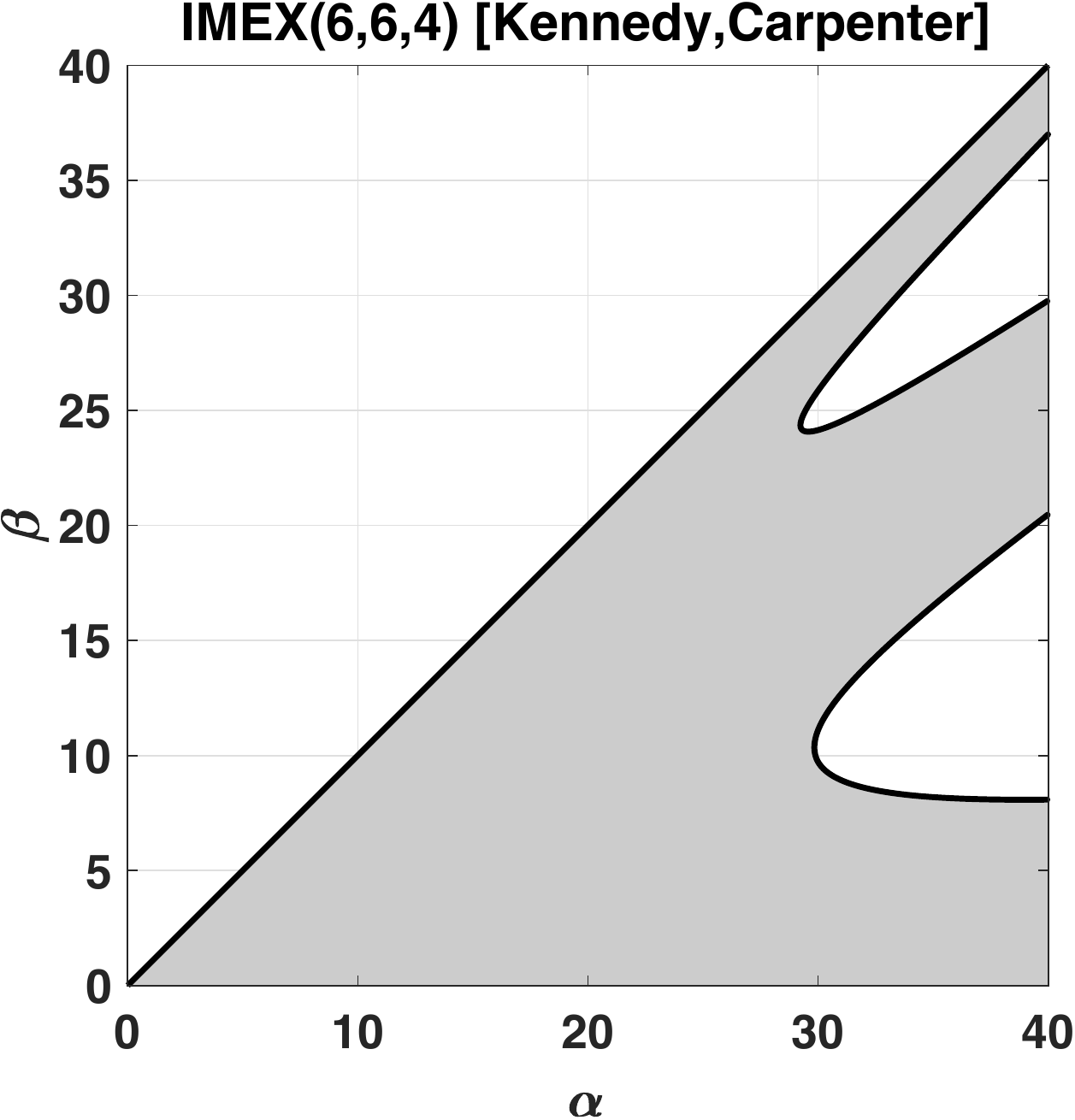}
 \caption{Stability regions of three \imex{} RK schemes for the scalar test
   problem $u_t = -\gamma u+ \mu u$. Left to right: 3rd order \imex{}(4,4,3);
           4th order \imex{} RK scheme by Cavaglieri and Bewley
           \cite{CavaglieriBewley2015};
           4th order \imex{} RK scheme by Kennedy and Carpenter
           \cite{KennedyCarpenter2003}.}
	\label{fig:ImExRK_StabRegion}
\end{figure}

A third-order \imex{} RK method that satisfies the wedge property
is the scheme denoted \imex{}(4,4,3) in \cite{AscherRuuthSpiteri1997}. It possesses 4 implicit stages and 4 explicit stages and has the Butcher tableau
\begin{equation}
\label{eq:imex443}
\begin{array}{c|c}
\myvec{c} & A \\
\hline
& \myvec{b}^T
\end{array}
=
\begin{array}{c|c c c c}
1/2 & 1/2 &   &  & \\
2/3 & 1/6 & 1/2 &  & \\
1/2 &-1/2 & 1/2 & 1/2 & \\
1   & 3/2 & -3/2& 1/2 & 1/2\\
\hline
& 3/2 & -3/2& 1/2 & 1/2
\end{array}
,\quad
\begin{array}{c|c}
\hat{\myvec{c}} & \hat{A} \\
\hline
& \hat{\myvec{b}}^T
\end{array}
=
\begin{array}{c|c c c c c}
0 & 0 &  &  &  & \\
1/2 & 1/2 & 0 &  &  & \\
2/3 & 11/18 & 1/18 & 0 &  & \\
1/2 & 5/6 & -5/6 & 1/2 & 0 & \\
1   & 1/4 & 7/4 & 3/4 & -7/4 & 0\\
\hline
& 1/4 & 7/4 & 3/4 & -7/4 & 0
\end{array}.
\end{equation}
Note that the explicit RK scheme $(\hat{A},\hat{\myvec{b}},\hat{\myvec{c}})$
is a 4-stage scheme recast as a 5-stage scheme which is stiffly accurate (the last row of
$\hat{A}$ equals the vector $\hat{\myvec{b}}^T$). This scheme is suitable for the PPE reformulation as it is unconditionally stable when solving the linear model problem \eqref{eq:testproblem} (see the left panel in Figure~\ref{fig:ImExRK_StabRegion} for the stability region). We therefore adopt the 3rd order \imex{} RK scheme \eqref{eq:imex443} for all numerical computations in this paper.

We are unaware of any existing 4th order (or higher) schemes that satisfy the wedge property. For instance, two popular 4th order \imex{} RK schemes, one by Cavaglieri and Bewley \cite{CavaglieriBewley2015} and the other by Kennedy and Carpenter \cite{KennedyCarpenter2003}, both violate the wedge property
(see Figure~\ref{fig:ImExRK_StabRegion} middle and right panels).

Note that the wedge property is only a sufficient condition for unconditional stability of the \imex{} RK scheme for the model problem \eqref{eq:testproblem}. Nevertheless, this condition can provide insight into the stability property of a given \imex{} RK scheme for the PPE reformulation (\ref{eq:SR-PPE_momentum}--\ref{eq:SR-PPE_pressure}). For instance, we generally observe that 4th order schemes violating the wedge property require a stiff parabolic time step restriction $\dt = O(\dx^2)$. 

As a final comment on the RK time discretization, RK schemes (including \imex{} RK schemes) suffer from order reduction when applied to initial boundary value problems with time-dependent boundary conditions or forcings on the boundary \cite{SanzVerwerHundsdorfer1986,Verwer1986,CarpenterGottliebAbarbanelDon1995}.  We summarize the order reduction phenomena in the following remark.

\begin{remark}
	\label{rem:OrderReduction}
 Order reduction is a generic problem of RK schemes where the observed temporal convergence rate is lower than the formal order of the scheme. Order reduction is often due to the formation of numerical boundary layers, caused by enforcing boundary conditions in the RK scheme \cite{RosalesSeiboldShirokoffZhou2017}. While several approaches exist to remedy order reduction 
 \cite{CarpenterGottliebAbarbanelDon1995,Alonso2002,RosalesSeiboldShirokoffZhou2017,KennedyCarpenter2003}, they do not easily generalize to PPE reformulations or IMEX settings. 
\end{remark}

To investigate the high order accuracy of the methods presented in this paper, we choose specific test problems (below) in which the structure of the problem does not lead to order reduction. Because order reduction is a generic phenomenon that is not specific to PPE reformulations, the (important) question of how to avoid order reduction is of generic nature and not specific to this work.

It is worth pointing out that \imex{} multistep methods are devoid of order reduction. A key obstacle to the applicability of multistep methods has been their restrictive stability properties, particularly for problems in which the implicit and explicit parts are both stiff \cite{SeiboldShirokoffZhou2018IMEXPractice}, such as PPE reformulations. However, recently proposed \imex{} multistep methods \cite{RosalesSeiboldShirokoffZhou2017IMEXTheory}
can achieve unconditional stability for such problems, and thus may provide an alternative suitable time-stepping strategy for PPE reformulations.

% --------------------------------------------------------------------------- 80
\subsection{Numerical Results for Manufactured Solutions}
\label{ssec:linear_NumerialResults}
% --------------------------------------------------------------------------- 80
The goal of this subsection is to demonstrate via numerical examples that, for the time-dependent Stokes problem (\ref{eq:linear_SR-PPE_u}--\ref{eq:linear_SR-PPE_p}), the proposed method can achieve high-order in space and third order in time. We employ the method of manufactured solutions and conduct convergence studies for the spatial accuracy (see \Srm\ref{sssec:linear_spatial}), the temporal accuracy for the 3rd order \imex{}(4,4,3) scheme \eqref{eq:imex443} (see \Srm\ref{sssec:linear_temp}), and the accuracy of the overall scheme (see \Srm\ref{sssec:linear_fullscheme}).

We implement the numerical methods using the software package \textsf{FEniCS} \cite{LoggMardalEtAl2012a}, which contains an extensive library of finite elements through one of its components: Finite element Automatic Tabulator (FIAT) \cite{Kirby2004a, Kirby2012a}. It provides $H(\tdiv)$ element spaces such as Raviart-Thomas (RT) elements and $H(\tcurl)$ elements of the N\'{e}d\'{e}lec types. The triangular meshes used for the computations are generated by the software package \textsf{Gmsh} \cite{GeuzaineRemacle2009}.

% --------------------------------------------------------------------------- 80
\subsubsection{Numerical Results for the Vector Heat Equation with EBC}
\label{sssec:vhe_results}
Before presenting the results for the PPE reformulation, we show the convergence results for the vector heat equation (VHE) with EBC as a benchmark. This provides both: a baseline for the PPE reformulation convergence study, and a verification of the code. Let the problem domain be $\Omega = [0,1)\times[0,1]$ with periodic b.c.\ applied in the $x$-direction and EBC in the $y$-direction. Hence, the problem reads as
\begin{subequations}
	\label{eq:vhe}
\begin{alignat}{2}
\myvec{u}_t &= \Delta\myvec{u} + \myvec{f} \quad &&\mbox{for~}(x,y)\in(0,1)^2,\\
\myvec{n}\times\myvec{u} &= 0,\;\;\DIV\myvec{u}=0 \quad &&\mbox{for~}(x,y)\in[0,1)\times\{0,1\},\\
\myvec{u}(0,y) &= \myvec{u}(1,y) \quad &&\mbox{for~}0<y<1.
\end{alignat}
\end{subequations}
The divergence-free manufactured solution $\myvec{u} = (u,v)^T = (\psi_y,-\psi_x)^T$ is generated by the stream function
\begin{equation*}
	\psi(x,y,t) = \cos(t)\sin(4\pi(x+y))(4y(1-y))^{4}.
\end{equation*}
The manufactured solution is selected so that the velocity field and its derivatives up to certain order vanish at the boundary to suppress the effect of order reduction due to the \imex{} RK time-stepping (Remark~\ref{rem:OrderReduction}). 

To conduct a spatial convergence study we use a series of regular meshes
with total number of elements 64, 256, 1024, 4096, and 16384. For the time
evolution we use the third order \imex{} (4,4,3) scheme with a fixed small time
step $\dt = 10^{-5}$ (so that the error due to the time-stepping is negligible)
and a fixed final time $T = 10^{-3}$ (100 time steps).
The FE spaces are chosen to be $(\sigma_h,\myvec{u}_h)\in P_r\times RT_{r-1}$,
with degrees $r = 1,2,3,4,5$. The observed spatial rates of convergence in
the $L^2$ norm are shown in
Table~\ref{tab:spatial_convergence_vhe} (the same rates are observed in the $L^{\infty}$ norm).
\begin{table}[htb!]
{\tiny
	\begin{tabular}{ccccccccc}
		\toprule
		\multicolumn{9}{l}{Spatial approximation order $r = 1$, $(\sigma_h,\myvec{u}_h)\in P_1\times RT_0$}\\
		\midrule
		$\dx$    & $\myvec{u}$ & Rate & $\DIV\myvec{u}$ & Rate &  $\sigma$ & Rate &  $\CURL\sigma$ & Rate\\ 
		\midrule
		2.50E-01 & 6.93E+00 & ---   & 9.97E+00 & ---  &  1.37E+02 & ---  &  2.76E+03 & ---  \\
		1.25E-01 & 3.26E+00 & 1.09  & 5.64E+00 & 0.82 &  3.19E+01 & 2.10 &  1.45E+03 & 0.92 \\
		6.25E-02 & 1.65E+00 & 0.98  & 1.04E+00 & 2.43 &  7.34E+00 & 2.12 &  7.24E+02 & 1.00 \\
		3.13E-02 & 8.24E-01 & 1.00  & 1.97E-01 & 2.41 &  1.77E+00 & 2.05 &  3.62E+02 & 1.00 \\
		1.56E-02 & 4.12E-01 & 1.00  & 4.57E-02 & 2.10 &  4.37E-01 & 2.01 &  1.81E+02 & 1.00 \\
		7.81E-03 & 2.06E-01 & 1.00  & 1.12E-02 & 2.03 &  1.09E-01 & 2.00 &  9.05E+01 & 1.00 \\
		\midrule
		\multicolumn{9}{l}{Spatial approximation order $r = 2$, $(\sigma_h,\myvec{u}_h)\in P_2\times RT_1$}\\
%		Spatial approximation order $r=2$ & & & & & & & & \\
		\midrule
		$\dx$    & $\myvec{u}$ & Rate & $\DIV\myvec{u}$ & Rate &  $\sigma$ & Rate &  $\CURL\sigma$ & Rate\\ 
		\midrule
		2.50E-01 & 2.13E+00 & ---  &  8.82E+00 & ---  & 2.46E+01 & ---  &  1.25E+03 & ---  \\
		1.25E-01 & 7.89E-01 & 1.43 &  1.27E+00 & 2.80 & 7.50E+00 & 1.71 &  3.97E+02 & 1.66 \\
		6.25E-02 & 2.13E-01 & 1.89 &  2.31E-01 & 2.46 & 1.99E+00 & 1.91 &  1.04E+02 & 1.94 \\
		3.13E-02 & 5.45E-02 & 1.97 &  3.09E-02 & 2.90 & 5.01E-01 & 1.99 &  2.62E+01 & 1.98 \\
		1.56E-02 & 1.37E-02 & 1.99 &  3.92E-03 & 2.98 & 1.25E-01 & 2.00 &  6.57E+00 & 2.00 \\
		\midrule
		\multicolumn{9}{l}{Spatial approximation order $r = 3$, $(\sigma_h,\myvec{u}_h)\in P_3\times RT_2$}\\
%		Spatial approximation order $r=3$ & & & & & & & & \\
		\midrule
		$\dx$    & $\myvec{u}$ & Rate & $\DIV\myvec{u}$ & Rate &  $\sigma$ & Rate &  $\CURL\sigma$ & Rate\\ 
		\midrule
		2.50E-01 & 1.32E+00 & ---  &  3.58E+00 & ---  &  1.88E+01 & ---  &  6.28E+02 & ---  \\
		1.25E-01 & 1.28E-01 & 3.37 &  2.27E-01 & 3.98 &  1.26E+00 & 3.90 &  7.40E+01 & 3.09 \\
		6.25E-02 & 1.56E-02 & 3.04 &  1.55E-02 & 3.87 &  7.96E-02 & 3.99 &  8.46E+00 & 3.13 \\
		3.13E-02 & 1.92E-03 & 3.02 &  9.67E-04 & 4.00 &  5.05E-03 & 3.98 &  1.05E+00 & 3.02 \\
		1.56E-02 & 2.40E-04 & 3.00 &  6.03E-05 & 4.00 &  3.17E-04 & 3.99 &  1.30E-01 & 3.00 \\
		\midrule
		\multicolumn{9}{l}{Spatial approximation order $r = 4$, $(\sigma_h,\myvec{u}_h)\in P_4\times RT_3$}\\
%		Spatial approximation order $r=4$ & & & & & & & & \\
		\midrule 
		$\dx$    & $\myvec{u}$ & Rate & $\DIV\myvec{u}$ & Rate &  $\sigma$ & Rate &  $\CURL\sigma$ & Rate\\ 
		\midrule
		2.50E-01 & 2.41E-01 & ---  & 1.68E+00 & ---  & 2.64E+00 & ---  & 1.52E+02 & ---   \\
		1.25E-01 & 1.74E-02 & 3.80 & 6.54E-02 & 4.68 & 1.70E-01 & 3.96 & 1.01E+01 & 3.91  \\
		6.25E-02 & 1.20E-03 & 3.85 & 1.04E-03 & 5.97 & 1.35E-02 & 3.66 & 6.98E-01 & 3.86  \\
		3.13E-02 & 7.73E-05 & 3.96 & 4.34E-05 & 4.58 & 8.98E-04 & 3.91 & 4.50E-02 & 3.95  \\
		1.56E-02 & 4.87E-06 & 3.99 & 2.27E-06 & 4.26 & 5.70E-05 & 3.98 & 2.84E-03 & 3.99  \\
		\midrule
		\multicolumn{9}{l}{Spatial approximation order $r = 5$, $(\sigma_h,\myvec{u}_h)\in P_5\times RT_4$}\\
%		Spatial approximation order $r=5$ & & & & & & & & \\
		\midrule
		$\dx$    & $\myvec{u}$ & Rate & $\DIV\myvec{u}$ & Rate &  $\sigma$ & Rate &  $\CURL\sigma$ & Rate\\ 
		\midrule
		2.50E-01 & 8.56E-02 & ---  &  4.14E-01 & ---  &  1.27E+00 & ---  &  4.39E+01 & ---   \\ 
		1.25E-01 & 2.16E-03 & 5.31 &  7.33E-03 & 5.82 &  2.86E-02 & 5.48 &  1.47E+00 & 4.90  \\ 
		6.25E-02 & 6.31E-05 & 5.10 &  5.62E-05 & 7.03 &  4.64E-04 & 5.94 &  4.13E-02 & 5.15  \\ 
		3.13E-02 & 1.92E-06 & 5.04 &  8.70E-07 & 6.01 &  7.52E-06 & 5.95 &  1.27E-03 & 5.03  \\
		1.56E-02 & 5.97E-08 & 5.01 &  1.36E-08 & 6.00 &  1.19E-07 & 5.99 &  3.93E-05 & 5.01  \\
		\bottomrule
	\end{tabular}
	}
\vspace{0.5em}
\caption{
	Observed spatial convergence rates in the $L^2$ norm for the vector heat equation
with electric boundary conditions \eqref{eq:vhe}. The spatial approximation orders are
$r=1,2,3,4,5$.}
\label{tab:spatial_convergence_vhe}
\end{table}
These results confirm the error estimate for the semi-discrete (in space) VHE in \cite{ArnoldChen2017} where the quantities $\myvec{u}$, $\DIV\myvec{u}$, $\sigma$ and $\CURL\sigma$ were proved to be at least $r$-th order convergent for $r$-th order finite elements.

The observed convergence rates are clean, with the exception of $r=4$ which shows some degradation in the convergence of $\nabla\cdot\myvec{u}$. Even though they do not contradict the existing error estimates, there are some convergence patterns that are worth commenting on:
\begin{enumerate}
\item The spatial convergence result for $\sigma$ exhibits an even-odd behavior. Specifically, the rate is $r$ for even order $r$, and it is $r+1$ when $r$ is odd. 
\item Extra orders of convergence (in space) for $\DIV\myvec{u}$ are observed for both the VPE and the VHE problems when the exact solution satisfies $\DIV\myvec{u}=0$. When the exact solutions are not divergence-free, there are no extra convergence orders for the divergence; however, $\DIV\myvec{u}=0$ is the common situation for incompressible fluid flows, so the extra order is noteworthy.
\end{enumerate}
Now that we have established the convergence results for the VHE with EBC, and compared them with the existing error estimates, we move on to the PPE reformulation. Note that, in this case, theoretical convergence results are not available. We numerically investigate the performance of the proposed schemes.

% --------------------------------------------------------------------------- 80
\subsubsection{Spatial Accuracy of the Time Dependent Stokes Problem}
\label{sssec:linear_spatial}
% --------------------------------------------------------------------------- 80
We now conduct a spatial convergence study of the proposed method for the time-dependent Stokes PPE (\ref{eq:linear_SR-PPE_u}--\ref{eq:linear_SR-PPE_p}) discretized in space according to \eqref{eq:mixed_SR-NSE_semi} and \eqref{eq:pressure_weak1} and in time via the \imex{} RK scheme \eqref{eq:imex443}.
We consider the same domain, $\Omega=[0,1)\times[0,1]$, with periodic b.c.\ in the $x$-direction, EBC in the $y$-direction, and the same divergence-free velocity profile, $\myvec{u}(x,y,t)$, as in \S\ref{sssec:vhe_results}. The pressure is taken to be
\begin{equation*}
	p = \cos(t)\cos(4\pi(x+y))(4y(1-y))^{4}\/.
\end{equation*}
The forcing is fixed as
$\myvec{f} = \myvec{u}_t -\nu\Delta\myvec{u}+\nabla p$, and the initial
conditions are chosen as $\myvec{u}_0 = \myvec{u}(x,y,0)$ (to match the
manufactured solution). Both the manufactured solution and the
forcing vanish at the boundary.

Following the same procedure as with the VHE above, we select the
FE space $(\sigma_h,\myvec{u}_h,p_h)\in P_r\times RT_{r-1}\times P_r$, with
degrees $r = 1,2,3,4,5$ on regular meshes. The time-stepping is done via
\imex{}(4,4,3), with fixed time step $\dt = 10^{-5}$ and final time $T = 10^{-3}$.
The stabilization parameter $\lambda$ is set to be 10, following the
suggestion in \cite{ShirokoffRosales2010}.

To avoid the approximation error introduced by iterative linear solvers, we use the sparse direct solver built into \textsf{FEniCS}. For the finest mesh resolution, and the 5th order spatial approximation, the sizes of the matrices for $(\sigma,\myvec{u})$ and $p$ are $656000 \times 656000$ and $205121 \times 205121$, respectively, which is close to the maximum that the direct solver could handle reliably.

\begin{table}[htb!]
{\tiny
	\begin{tabular}{ccccccccccc}
		\toprule
		\multicolumn{11}{l}{Spatial approximation order $r = 1$, $(\sigma_h,\myvec{u}_h,p_h)\in P_1\times RT_0\times P_1$}\\
		\midrule
		$\dx$    & $\myvec{u}$ & Rate & $\DIV\myvec{u}$ & Rate & $\nabla\myvec{u}$ & Rate & $\sigma$ & Rate & $p$ & Rate  \\
		\midrule
		2.50E-01 & 6.93E+00 & ---- & 1.45E+01 & ---- & 1.37E+02 & ---- & 1.37E+02 & ---- & 3.39E+01 & ----  \\ 
		1.25E-01 & 3.26E+00 & 1.09 & 6.13E+00 & 1.25 & 1.37E+02 & 0.00 & 3.19E+01 & 2.10 & 3.38E+00 & 3.33  \\ 
		6.25E-02 & 1.65E+00 & 0.98 & 1.05E+00 & 2.54 & 1.37E+02 & 0.00 & 7.33E+00 & 2.12 & 4.41E-01 & 2.94  \\ 
		3.13E-02 & 8.24E-01 & 1.00 & 1.96E-01 & 2.42 & 1.37E+02 & 0.00 & 1.77E+00 & 2.05 & 1.08E-01 & 2.03  \\
		1.56E-02 & 4.12E-01 & 1.00 & 4.57E-02 & 2.10 & 1.37E+02 & 0.00 & 4.37E-01 & 2.01 & 2.71E-02 & 1.99  \\
		7.81E-03 & 2.06E-01 & 1.00 & 1.12E-02 & 2.02 & 1.37E+02 & 0.00 & 1.09E-01 & 2.00 & 6.78E-03 & 2.00  \\
		\midrule
		\multicolumn{11}{l}{Spatial approximation order $r = 2$, $(\sigma_h,\myvec{u}_h,p_h)\in P_2\times RT_1\times P_2$}\\
		\midrule
		$\dx$    & $\myvec{u}$ & Rate & $\DIV\myvec{u}$ & Rate & $\nabla\myvec{u}$ & Rate & $\sigma$ & Rate & $p$ & Rate  \\
		\midrule
		2.50E-01 & 2.09E+00 & ---- & 7.40E+00 & ---- & 9.33E+01 & ---- & 2.46E+01 & ----  & 1.80E+01 & ----  \\ 
		1.25E-01 & 7.88E-01 & 1.41 & 1.30E+00 & 2.51 & 6.25E+01 & 0.58 & 7.50E+00 & 1.71  & 7.64E-01 & 4.56  \\ 
		6.25E-02 & 2.13E-01 & 1.89 & 2.41E-01 & 2.43 & 3.27E+01 & 0.93 & 1.99E+00 & 1.91  & 6.18E-02 & 3.63  \\ 
		3.13E-02 & 5.45E-02 & 1.97 & 3.12E-02 & 2.95 & 1.66E+01 & 0.98 & 5.01E-01 & 1.99  & 1.05E-02 & 2.56  \\ 
		1.56E-02 & 1.37E-02 & 1.99 & 3.93E-03 & 2.99 & 8.33E+00 & 0.99 & 1.25E-01 & 2.00  & 2.45E-03 & 2.09  \\ 
		7.81E-03 & 3.43E-03 & 2.00 & 4.92E-04 & 3.00 & 4.17E+00 & 1.00 & 3.13E-02 & 2.00  & 6.08E-04 & 2.01  \\
		\midrule
		\multicolumn{11}{l}{Spatial approximation order $r = 3$, $(\sigma_h,\myvec{u}_h,p_h)\in P_3\times RT_2\times P_3$}\\
		\midrule
		$\dx$    & $\myvec{u}$ & Rate & $\DIV\myvec{u}$ & Rate & $\nabla\myvec{u}$ & Rate & $\sigma$ & Rate & $p$ & Rate  \\
		\midrule
		2.50E-01 & 1.33E+00 & ---- & 3.81E+00 & ---- & 7.09E+01 & ---- & 1.88E+01 & ----  & 4.33E+00 & ----  \\ 
		1.25E-01 & 1.28E-01 & 3.38 & 2.27E-01 & 4.07 & 1.67E+01 & 2.09 & 1.26E+00 & 3.90  & 1.46E-01 & 4.89  \\ 
		6.25E-02 & 1.56E-02 & 3.04 & 1.51E-02 & 3.92 & 4.24E+00 & 1.98 & 7.97E-02 & 3.99  & 6.91E-03 & 4.40  \\ 
		3.13E-02 & 1.92E-03 & 3.02 & 9.58E-04 & 3.97 & 1.06E+00 & 2.00 & 5.05E-03 & 3.98  & 2.96E-04 & 4.54  \\ 
		1.56E-02 & 2.40E-04 & 3.00 & 6.01E-05 & 4.00 & 2.66E-01 & 2.00 & 3.17E-04 & 4.00  & 1.47E-05 & 4.34  \\ 
		7.81E-03 & 2.99E-05 & 3.00 & 3.76E-06 & 4.00 & 6.65E-02 & 2.00 & 1.98E-05 & 4.00  & 9.81E-07 & 3.90  \\
		\midrule
		\multicolumn{11}{l}{Spatial approximation order $r = 4$, $(\sigma_h,\myvec{u}_h,p_h)\in P_4\times RT_3\times P_3$}\\
		\midrule
		$\dx$    & $\myvec{u}$ & Rate & $\DIV\myvec{u}$ & Rate & $\nabla\myvec{u}$ & Rate & $\sigma$ & Rate & $p$ & Rate  \\
		\midrule
		2.50E-01 & 2.41E-01 & ---- & 1.66E+00 & ---- & 1.94E+01 & ---- & 2.66E+00 & ---- & 1.17E+00 & ----  \\ 
		1.25E-01 & 1.74E-02 & 3.79 & 6.23E-02 & 4.74 & 3.03E+00 & 2.68 & 1.75E-01 & 3.92 & 4.83E-02 & 4.60  \\ 
		6.25E-02 & 1.20E-03 & 3.86 & 1.01E-03 & 5.95 & 4.04E-01 & 2.91 & 1.38E-02 & 3.67 & 3.00E-03 & 4.01  \\ 
		3.13E-02 & 7.73E-05 & 3.96 & 4.31E-05 & 4.55 & 5.14E-02 & 2.97 & 9.03E-04 & 3.93 & 1.12E-04 & 4.74  \\ 
		1.56E-02 & 4.87E-06 & 3.99 & 2.26E-06 & 4.25 & 6.46E-03 & 2.99 & 5.71E-05 & 3.98 & 4.62E-06 & 4.60  \\	
		\midrule
		\multicolumn{11}{l}{Spatial approximation order $r = 5$, $(\sigma_h,\myvec{u}_h,p_h)\in P_5\times RT_4\times P_5$}\\
		\midrule
		$\dx$    & $\myvec{u}$ & Rate & $\DIV\myvec{u}$ & Rate & $\nabla\myvec{u}$ & Rate & $\sigma$ & Rate & $p$ & Rate  \\
		\midrule
		2.50E-01& 8.61E-02 & --- & 4.24E-01 & --- & 6.56E+00 & --- & 1.31E+00 & --- & 4.88E-01 & --- \\
		1.25E-01& 2.17E-03 & 5.31& 7.49E-03 & 5.82& 4.80E-01 & 3.77& 3.10E-02 & 5.40& 1.36E-02 & 5.16\\
		6.25E-02& 6.33E-05 & 5.10& 5.48E-05 & 7.10& 3.14E-02 & 3.94& 7.41E-04 & 5.39& 5.86E-04 & 4.54\\
		3.13E-02& 1.93E-06 & 5.04& 8.73E-07 & 5.97& 1.98E-03 & 3.99& 2.31E-05 & 5.00& 2.19E-05 & 4.74\\
		1.56E-02& 5.97E-08 & 5.01& 1.53E-08 & 5.83& 1.24E-04 & 4.00& 6.81E-07 & 5.09& 6.70E-07 & 5.03\\
		\bottomrule
	\end{tabular}
}
\vspace{0.5em}
\caption{Spatial error convergence in the $L^2$ norm for the time dependent Stokes equation (\ref{eq:linear_SR-PPE_u}--\ref{eq:linear_SR-PPE_p}) in \S\ref{sssec:linear_spatial}. The spatial approximations are $(\sigma_h,\myvec{u}_h,p_h)\in P_r\times RT_{r-1}\times P_r$, with degrees $r = 1,2,3,4,5$. The calculations are done using $\lambda = 10$ and the 3rd order \imex{}(4,4,3), with a fixed small time step $\dt = 10^{-5}$ and a fixed final time $T = 10^{-3}$.}
	 \label{tab:spatial_err_cvg}
\end{table}

Table~\ref{tab:spatial_err_cvg} shows the spatial error convergence results for different degrees of the spatial approximation. The approximation errors are measured in the $L^2$ norm. For quantities related to $\myvec{u}$ and $\sigma$, the rates of convergence have similar behaviors as for the VHE. Specifically: (i)~the velocity $\myvec{u}$ is $r$-th order convergent for an $r$-th order spatial approximation, and (ii)~the extra convergence orders for $\DIV\myvec{u}$ is carried over to the PPE reformulation.
The convergence rate for the pressure $p$ behaves similarly to $\sigma$, which
appears in the right hand side of the weak formulation for the pressure in
\eqref{eq:pressure_weak1} and \eqref{eq:pressure_weak2}. However, non-clean
pressure convergence rates are observed for approximation orders larger then 2.

It should be stressed that we numerically measure the convergence rate for the error in the velocity gradient $\nabla\myvec{u}$, even though the FEM spaces do not guarantee that $\nabla\myvec{u}_h$ is in $L^2$. In particular, the numerical solution for the velocity is, generally, discontinuous (in the tangential direction) across Raviart-Thomas element edges. Here we measure the error in $\nabla\myvec{u}_h$ by ignoring the jumps across the edges, i.e.\, by only counting the error within each element.  Note that it is of interest to measure the accuracy of velocity gradients as they relate to fluid stresses and forces at the boundary of objects. In addition, measuring the accuracy of fluid gradients will be a precursor to \S\ref{sec:NonlinearAdvection}, in which $\myvec{u}\cdot\nabla \myvec{u}$ will be included in the equations.

% --------------------------------------------------------------------------- 80
\subsubsection{Temporal Accuracy of the Time Dependent Stokes Problem}
\label{sssec:linear_temp}
% --------------------------------------------------------------------------- 80
To check the temporal accuracy, we consider a manufactured solution
on the same domain $\Omega=[0,1)\times[0,1]$ with periodicity in the
$x$-direction. The manufactured solution follows from the stream function
\begin{equation*}
	\psi(x,y,t) = \cos(200t)\sin^2(\pi x)\sin^2(\pi y)\/,
\end{equation*}
so that the divergence-free condition is automatically satisfied by the
velocity field $\myvec{u} = (u,v)^T = (\psi_y, -\psi_x)^T$. The pressure is
\begin{equation*}
	p(x,y,t) = \cos(200t)\sin(2\pi x)\sin(\pi y)\/.
\end{equation*}
Again, note that both the solution and the forcing vanish at the domain
boundary. To test the temporal errors, we select a highly oscillatory in time manufactured solution. The high frequency oscillations in time ensure that the time discretization errors dominate the spatial discretization errors for the mesh resolution we use.

To perform the temporal error convergence study we solve the linear
problem on a fixed mesh (mesh size $\dx = 3.125\times 10^{-2}$ and
4096 elements), a fixed FE discretization scheme (FE space with
degree 4, that is $(\sigma_h,\myvec{u}_h,p_h)\in P_4\times RT_3\times P_4$),
and a final time $T = 0.5$.
Then we vary the time step: $\dt = 2^{-k}$, $k = 7,\dots,13$.

The temporal convergence results for the 3rd order \imex{}(4,4,3) scheme
are shown in Table~\ref{tab:temporal_err_cvg}. All quantities, except
for the divergence $\DIV\myvec{u}$, exhibit a 3rd order convergence in
time in the $L^2$ norm, while $\DIV\myvec{u}$ remains small for all
$\dt$. The stagnation in the convergence for $\nabla\myvec{u}$, at
$O(10^{-4})$, is due to
the spatial approximation error dominating the temporal error.
% =========================================================================== 80

\begin{table}[htb!]
{\tiny
	\begin{tabular}{ccccccccccc}
		\toprule
%		\midrule
		$\dt$   & $\myvec{u}$  & Rate & $\DIV\myvec{u}$ & Rate & $\nabla\myvec{u}$ & Rate & $\sigma$ & Rate & $p$ & Rate     \\ 
		\midrule
		7.81E-03& 7.725935E-02 & ---  & 6.531237E-08 &  ---  & 5.932904E-01 & ---   & 5.932850E-01 & ---  & 9.039549E-02 & ---  \\
		3.91E-03& 9.910506E-03 &  2.96  & 3.476197E-08 &  (0.91) & 7.625508E-02 & 2.96  & 7.624960E-02 &  2.96  & 2.074725E-02 & 2.12 \\
		1.95E-03& 1.269255E-03 &  2.97  & 3.004490E-08 &  (0.21) & 9.693504E-03 & 2.98  & 9.685865E-03 &  2.98  & 3.396972E-03 & 2.61 \\ 
		9.77E-04& 1.616344E-04 &  2.97  & 3.032122E-08 & (-0.01) & 1.233570E-03 & 2.97  & 1.207865E-03 &  3.00  & 4.883108E-04 & 2.80 \\ 
		4.88E-04& 2.073504E-05 &  2.96  & 3.101061E-08 & (-0.03) & 2.734583E-04 & (2.17)  & 1.495803E-04 &  3.01 & 6.545107E-05 & 2.90 \\
		2.44E-04& 2.760336E-06 &  2.91  & 3.146726E-08 & (-0.02) & 2.270653E-04 & (0.27)  & 1.894479E-05 &  2.98 & 8.380049E-06 & 2.97 \\ 
		1.22E-04& 4.870641E-07 & (2.50) & 3.171521E-08 & (-0.01) & 2.259969E-04 & (0.01)  & 2.585832E-06 & (2.87) & 1.032259E-06 & 3.02 \\ 
		\bottomrule
	\end{tabular}
}
\vspace{0.5em}
\caption{Temporal error convergence for (\ref{eq:linear_SR-PPE_u}--\ref{eq:linear_SR-PPE_p}) in the $L^2$ norm for the 3rd order \imex{}(4,4,3), on a fixed mesh, with the 4th order spatial discretization
	$P_4\times RT_3 \times P_4$. The convergence rates in parentheses stagnate due to the spatial error dominating the temporal error.}
\label{tab:temporal_err_cvg}
\end{table}

% --------------------------------------------------------------------------- 80
%
% --------------------------------------------------------------------------- 80
\subsubsection{Convergence Results in Both Space and Time for the Time Dependent Stokes Equations}
\label{sssec:linear_fullscheme}
% --------------------------------------------------------------------------- 80
Here we present the convergence results for a 3rd order scheme in both space and time applied to the time-dependent Stokes problem (\ref{eq:linear_SR-PPE_u}--\ref{eq:linear_SR-PPE_p}). We use the 3rd order spatial discretization $(\sigma_h,\myvec{u}_h,p_h)\in P_3\times RT_{2}\times P_3$ and the 3rd order time-stepping \imex{}(4,4,3), with the time step scaled proportional to the mesh size ---
specifically $\dt = 0.2\dx$. A series of regular meshes with total number of elements 64, 256, 1024, 4096, 16384 and 65536 are used. On the finest mesh, the sizes of the linear systems are 983808 for $(\myvec{u},\sigma)$ and 295297 for $p$.
The manufactured solution is chosen to be
\begin{align*}
	u = \psi_y,\quad
	v = -\psi_x,\quad
	p = \pi\cos(t)\cos(\pi x)\sin(\pi y)(4x(1-x))^4(4y(1-y))^4, 
\end{align*}
where $\psi(x,y,t) = \cos(t)\sin^2(\pi x)\sin^2(\pi y)(4x(1-x))^4(4y(1-y))^4$,
on the unit square domain $\Omega=[0,1]^2$. The EBC are prescribed at all boundaries. \emph{This test differs from the previous tests where periodic b.c.\ in $x\/$ were used.} Again, the exact solution is constructed so that the velocity and the forcing vanish at the boundary. 

Table~\ref{tab:err_cvg_linear} shows the error convergence results in the $L^2$ norm (top) and in the $L^{\infty}$ norm (bottom). Clean convergence
are observed for $\myvec{u}$, $\nabla\myvec{u}$ and $\DIV\myvec{u}$, with rates 3, 2 and 4 respectively. Non-clean convergence rates occur for $\sigma$, $\CURL\sigma$, $p$ and $\nabla p$. However, they appear to be close to 3.  

\begin{table}[htb!]
	{\tiny
		\begin{tabular}{ccccccccccccccc}
			\toprule
			\multicolumn{15}{l}{Error convergence in the $L^2$ norm}\\
			\midrule
			$\dx$   & $\myvec{u}$ & Rate & $\DIV\myvec{u}$ & Rate & $\nabla\myvec{u}$ & Rate & $\sigma$ & Rate & $\CURL\sigma$ & Rate & $p$ & Rate & $\nabla p$   & Rate \\ 
			\midrule
			2.50E-01& 4.27E-02 & --- & 9.84E-02 & --- & 2.91E+00 & --- & 3.22E-01 & --- & 1.40E+01 & --- & 1.99E-01 & --- & 5.37E+00 & --- \\
			1.25E-01& 8.12E-03 & 2.39& 8.25E-03 & 3.57& 1.03E+00 & 1.51& 2.48E-02 & 3.70& 2.08E+00 & 2.75& 2.76E-03 & 6.17& 2.11E-01 & 4.67\\
			6.25E-02& 1.04E-03 & 2.96& 5.87E-04 & 3.81& 2.62E-01 & 1.97& 2.81E-03 & 3.14& 2.76E-01 & 2.91& 3.53E-04 & 2.96& 3.09E-02 & 2.77\\
			3.13E-02& 1.33E-04 & 2.97& 3.76E-05 & 3.96& 6.57E-02 & 2.00& 4.57E-04 & 2.62& 3.56E-02 & 2.96& 9.88E-05 & 1.84& 5.14E-03 & 2.59\\
			1.56E-02& 1.69E-05 & 2.97& 2.37E-06 & 3.99& 1.64E-02 & 2.00& 7.66E-05 & 2.58& 4.56E-03 & 2.96& 2.33E-05 & 2.08& 6.73E-04 & 2.93\\
			7.81E-03& 2.14E-06 & 2.99& 1.49E-07 & 3.99& 4.11E-03 & 2.00& 1.14E-05 & 2.74& 5.83E-04 & 2.97& 2.91E-06 & 3.00& 7.65E-05 & 3.13\\
			\midrule
			\multicolumn{15}{l}{Error convergence in the $L^\infty$ norm}\\
			\midrule
			$\dx$   & $\myvec{u}$ & Rate & $\DIV\myvec{u}$ & Rate & $\nabla\myvec{u}$ & Rate & $\sigma$ & Rate & $\CURL\sigma$ & Rate & $p$ & Rate & $\nabla p$   & Rate \\ 
			\midrule
			2.50E-01& 2.42E-01 & ---& 3.10E-01 & ---& 1.30E+01 & ---& 2.32E+00 & ---& 1.05E+02 & ---& 1.23E+00 & ---& 2.87E+01 & --- \\
			1.25E-01& 1.12E-01 & 1.12& 4.46E-02 & 2.80& 7.80E+00 & 0.74& 1.54E-01 & 3.91& 2.26E+01 & 2.21& 3.44E-02 & 5.16& 1.45E+00 & 4.31 \\
			6.25E-02& 1.77E-02 & 2.65& 3.17E-03 & 3.81& 2.12E+00 & 1.88& 1.92E-02 & 3.01& 2.52E+00 & 3.16& 3.99E-03 & 3.11& 4.28E-01 & 1.76 \\
			3.13E-02& 2.34E-03 & 2.92& 2.06E-04 & 3.95& 5.36E-01 & 1.98& 2.68E-03 & 2.84& 3.66E-01 & 2.79& 7.06E-04 & 2.50& 9.46E-02 & 2.18 \\
			1.56E-02& 2.97E-04 & 2.98& 1.29E-05 & 4.00& 1.34E-01 & 2.00& 4.05E-04 & 2.73& 5.34E-02 & 2.78& 1.26E-04 & 2.48& 1.42E-02 & 2.74 \\
			7.81E-03& 3.73E-05 & 2.99& 8.00E-07 & 4.01& 3.35E-02 & 2.00& 6.02E-05 & 2.75& 5.91E-03 & 3.18& 1.52E-05 & 3.05& 2.29E-03 & 2.63 \\
			\bottomrule
		\end{tabular}
	}
\vspace{0.5em}
\caption{Error convergence for the linear problem without the advection
	term  (\ref{eq:linear_SR-PPE_u}--\ref{eq:linear_SR-PPE_p}) in the $L^2$ norm (top) and the $L^{\infty}$ norm (bottom). The problem
	is solved with the 3rd order \imex{} RK scheme (\imex{}(4,4,3)), the 3rd order
	spatial discretization $P_3\times RT_2\times P_3$, and $\lambda=30$.}
\label{tab:err_cvg_linear}
\end{table}

\begin{figure}[htb!]
	\includegraphics[width=0.45\textwidth]
	{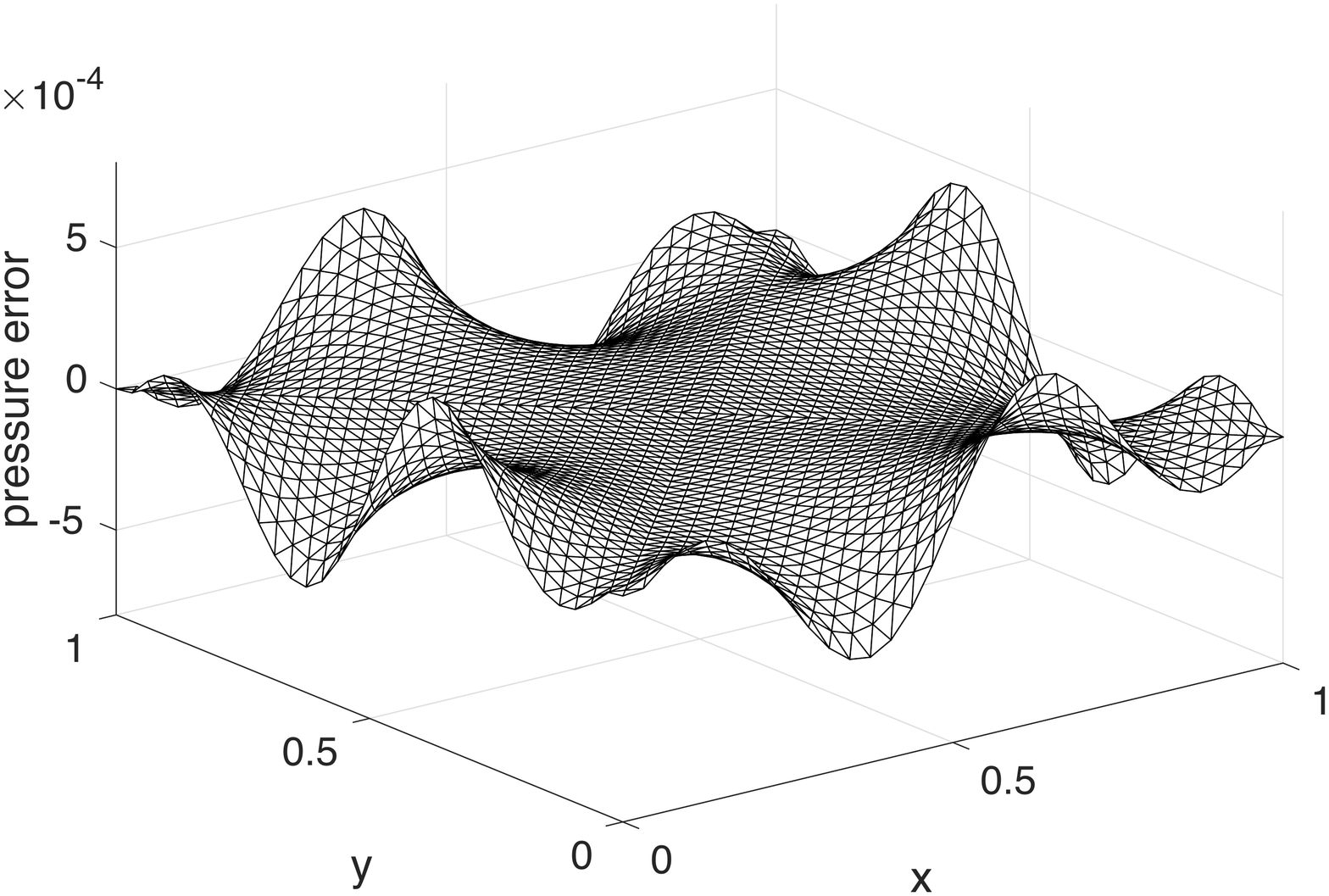}\quad
	\includegraphics[width=0.45\textwidth]
	{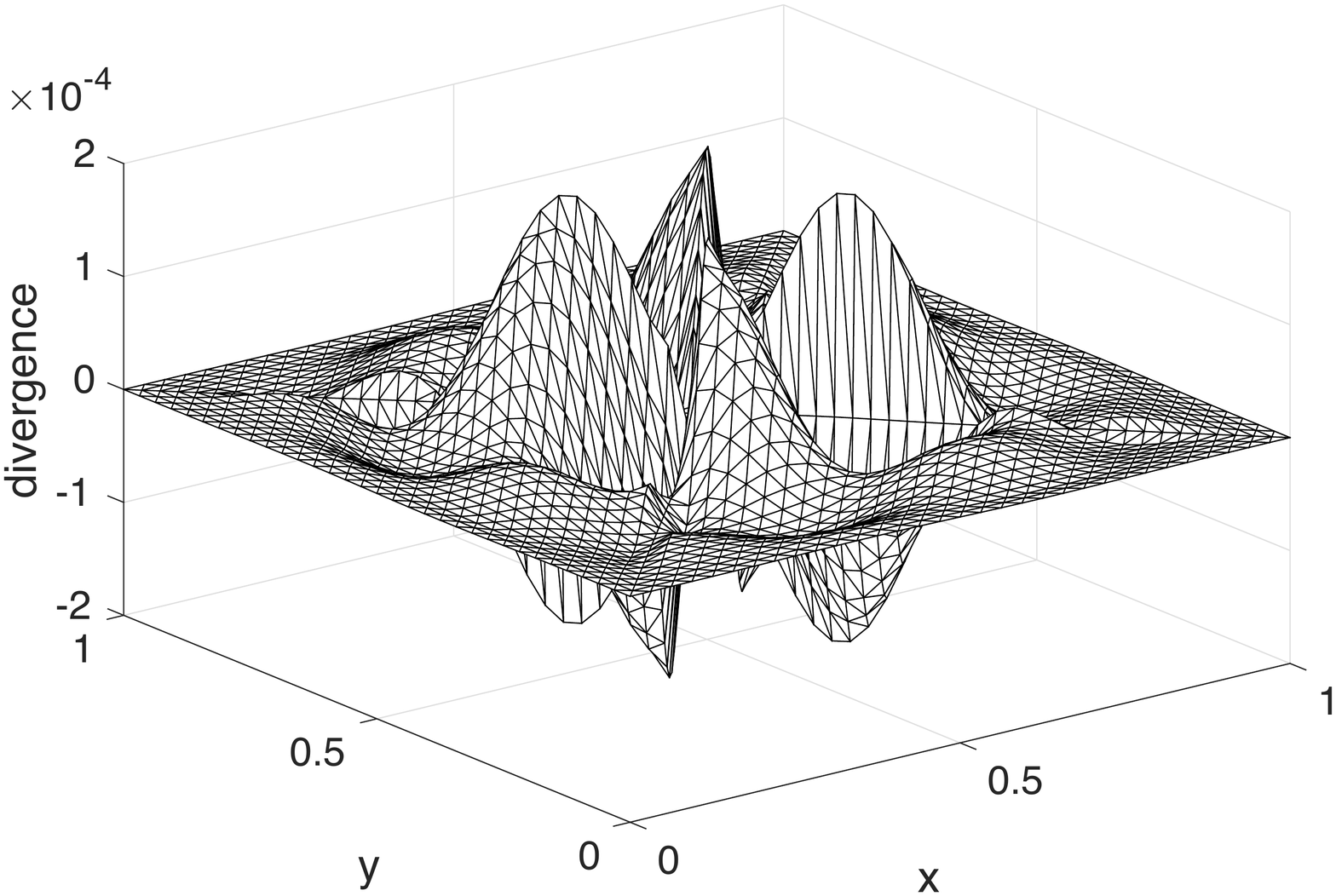}
\caption{Errors in the pressure (left) and the divergence (right) for the time dependent Stokes equations in \S\ref{sssec:linear_fullscheme}, as functions of $(x,y)$, at the final time $T=3$. 
The plot is for the 3rd order \imex{}(4,4,3) with a 3rd order ($r = 3$)
 spatial discretization on a regular mesh with 4096 elements.}
 \label{fig:pressure_div_error}
\end{figure}

Figure~\ref{fig:pressure_div_error} shows the pressure error and the divergence, as functions of
$(x,y)$, at the final time $T=3$. No numerical boundary layers are observed in the pressure error.
The divergence at the final time is small at the boundary, but not exactly zero. This is due to the fact
that the divergence boundary condition is only enforced weakly in the mixed
formulation (appears as a natural boundary condition). The dynamics of the
PPE reformulation ($\DIV\myvec{u}$ satisfies a heat equation) keeps the
divergence small across the domain.

In the PPE reformulation, there is no inf-sup condition for the velocity and
the pressure. This is in contrast to the conventional FEM formulation for the Navier-Stokes
equations, in which the finite element spaces for the velocity and the
pressure need to satisfy the discrete version of the inf-sup condition to ensure stability.
Hence the orders of the spatial approximations for
$(\sigma_h,\myvec{u}_h)\in P_r\times RT_{r-1}$ and the pressure $p_h\in P_r$
could be selected independently. However, we observed in the numerical
experiments that: \emph{if the pressure approximation is one order less than
the velocity approximation, the error convergence for the divergence
decreases by one compared to the equal-order case.} Furthermore:
\emph{increasing the order of approximation for the pressure does not
 improve the convergence rates compared to the equal-order case}. 

\begin{remark}(Degradation of convergence rates)
	Some of the convergence rates in Table~\ref{tab:err_cvg_linear} do not
	exhibit a clear integer rate.  A possible explanation for the
	reduction in convergence rate is provided in
	\cite{ArnoldFalkGopalakrishnan2012} (Theorems 5.1--5.2),
	where it is shown that a mixed finite element approximation
        for the Stokes equations exhibits a degraded convergence rate.
\end{remark}
% --------------------------------------------------------------------------- 80
%
% --------------------------------------------------------------------------- 80
\subsubsection{Influence of the parameter $\lambda$}
\label{sssec:param_lambda}
% --------------------------------------------------------------------------- 80
The role of the parameter $\lambda$ in the PPE reformulation (\ref{eq:SR-PPE_momentum}--\ref{eq:SR-PPE_pressure}) is to
exponentially enforce the normal velocity condition at the boundary, as highlighted by the ODE \eqref{eq:SR-normal}.
If $\lambda=0$, the exact solution to the ODE \eqref{eq:SR-normal} still satisfies $\myvec{n}\cdot\myvec{u}=\myvec{n}\cdot\myvec{g}$ (when $\phi = 0$), however the numerical scheme exhibits a steady growth of the error, as described in \cite{ShirokoffRosales2010}. Our numerical experiments show that when $\lambda$ is set to 0, there are
degradations in the convergence rates for $\myvec{u}$, $\sigma$ and the
pressure $p$, even though a 3rd order scheme is used. Increasing $\lambda$
controls the error in the velocity and hence improves the accuracy of the
overall scheme (see Table~\ref{tab:err_cvg_lambda}).

\begin{table}[htb!]
	{\tiny
		\begin{tabular}{ccccccccccc}
			\toprule
			\multicolumn{11}{l}{Error convergence in the $L^2$ norm with $\lambda = 0$}\\
			\midrule
			$\dx$   & $\myvec{u}$ & Rate & $\DIV\myvec{u}$ & Rate & $\nabla\myvec{u}$ & Rate & $\sigma$ & Rate & $p$ & Rate\\ 
			\midrule
			2.50E-01& 5.61E-02 & --- & 9.60E-02 & --- & 2.95E+00 & --- & 5.06E-01 & --- & 4.39E-01 & --- \\
			1.25E-01& 3.47E-02 & 0.69& 8.24E-03 & 3.54& 1.06E+00 & 1.48& 2.51E-01 & 1.01& 1.88E-01 & 1.23\\
			6.25E-02& 1.81E-02 & 0.94& 5.86E-04 & 3.81& 2.97E-01 & 1.84& 1.35E-01 & 0.89& 1.07E-01 & 0.81\\
			3.13E-02& 5.65E-03 & 1.68& 3.76E-05 & 3.96& 7.99E-02 & 1.89& 4.51E-02 & 1.59& 3.86E-02 & 1.47\\
			1.56E-02& 9.61E-04 & 2.55& 2.37E-06 & 3.99& 1.89E-02 & 2.08& 9.26E-03 & 2.28& 8.58E-03 & 2.17\\			
			\midrule
			\multicolumn{11}{l}{Error convergence in the $L^2$ norm with $\lambda = 30$}\\
			\midrule
			$\dx$   & $\myvec{u}$ & Rate & $\DIV\myvec{u}$ & Rate & $\nabla\myvec{u}$ & Rate & $\sigma$ & Rate & $p$ & Rate\\ 
			\midrule
			2.50E-01& 4.27E-02 & --- & 9.84E-02 & --- & 2.91E+00 & --- & 3.22E-01 & --- & 1.99E-01 & --- \\
			1.25E-01& 8.12E-03 & 2.39& 8.25E-03 & 3.57& 1.03E+00 & 1.51& 2.48E-02 & 3.70& 2.76E-03 & 6.17\\
			6.25E-02& 1.04E-03 & 2.96& 5.87E-04 & 3.81& 2.62E-01 & 1.97& 2.81E-03 & 3.14& 3.53E-04 & 2.96\\
			3.13E-02& 1.33E-04 & 2.97& 3.76E-05 & 3.96& 6.57E-02 & 2.00& 4.57E-04 & 2.62& 9.88E-05 & 1.84\\
			1.56E-02& 1.69E-05 & 2.97& 2.37E-06 & 3.99& 1.64E-02 & 2.00& 7.66E-05 & 2.58& 2.33E-05 & 2.08\\
			7.81E-03& 2.14E-06 & 2.99& 1.49E-07 & 3.99& 4.11E-03 & 2.00& 1.14E-05 & 2.74& 2.91E-06 & 3.00\\
			\bottomrule
		\end{tabular}
	}
	\vspace{0.5em}
	\caption{Error convergence in the $L^2$ norm for the linear problem 
	(\ref{eq:linear_SR-PPE_u}--\ref{eq:linear_SR-PPE_p}) with $\lambda=0$ (top)
		and $\lambda=30$ (bottom). The problem is solved with the 3rd order
		approximation $P_3\times RT_2\times P_3$ and the 3rd order \imex{}(4,4,3),
		using the same manufactured solution in \Srm\ref{sssec:linear_fullscheme}.}
	\label{tab:err_cvg_lambda}
\end{table}

It is worth mentioning that the convergence rate of $\DIV\myvec{u}$, even
the error in $\DIV\myvec{u}$, seems to be unaffected by the choice of
$\lambda$.
This is due to the fact that the enforcement of the divergence condition
derives from the fact that the divergence (implicitly) satisfies the heat
equation \eqref{eq:SR-heat}, with homogeneous boundary condition. Hence it
is independent of the parameter $\lambda$.
When $\lambda = 0\/$ the normal velocity
condition is enforced poorly, thus there is (numerical) flow through the
boundary. But the divergence is small, hence the total flow through the
boundary should be small, with inflows compensated by outflows.
%
% --------------------------------------------------------------------------- 80
\section{Treatment of the Nonlinear Advection Term}
\label{sec:NonlinearAdvection}
% --------------------------------------------------------------------------- 80

Next, in \Srm\ref{ssec:NonlinearAdvection_discretization}, we explore avenues
for discretizing the nonlinear advection term in the mixed finite element
framework introduced in \Srm\ref{sec:method_linear}. Although the approach
we present is ad hoc, the numerical results in
\Srm\ref{ssec:NonlinearAdvection_results} indicate that the full scheme is
convergent, with a little degradation in the rates of convergence relative to
the linearized equations. Benchmark test results for the lid-driven cavity
and the backward-facing step flow are presented in
\Srm\ref{ssec:NonlinearAdvection_benchmark}, which show good agreement with
the reference data.
% --------------------------------------------------------------------------- 80
%
% --------------------------------------------------------------------------- 80
\subsection{Discretization of the nonlinear advection term}
\label{ssec:NonlinearAdvection_discretization}
% --------------------------------------------------------------------------- 80
In \S\ref{sec:method_linear}, we introduced a mixed finite element spatial discretization as a way to 
handle the electric boundary conditions (EBC) for the vector heat equation and time-dependent Stokes equation. While the mixed finite elements resolve several difficulties for the EBC, they come with a caveat:  $RT$ elements approximating $\myvec{u}$ are only guaranteed to be continuous across interior edges in the normal direction but can jump in the tangential direction. Therefore, representations of the discrete solution using $RT$ elements are only weakly differentiable along the normal direction across edges.  This creates a problem then for handling nonlinear advection terms $\NL{u}{}=\nl{\myvec{u}}{}{}$, and in fact, even linear advection terms $(\myvec{a}\cdot\nabla)\myvec{u}$. Specifically, in 2D, $RT$ elements approximate a function $\myvec{u}$ $\myvec{v}=(v_1,v_2)^T\in H(\tdiv;\Omega)$, i.e.\, $v_1, v_2, \DIV\myvec{v}\in L^2(\Omega)$.  In general, each component of
$\nabla\myvec{v}$ is not guaranteed to be in $L^2$, therefore for any
$\myvec{u},\myvec{v}\in H(\tdiv;\Omega)$ the inner product
\begin{equation}
\label{eq:adv_integral}
\inprod{(\myvec{a}\cdot\nabla)\myvec{u}}{\myvec{v}} =
 \int_{\Omega}(\myvec{a}\cdot\nabla)\myvec{u}\cdot\myvec{v}\intd V
\end{equation}
is not properly defined.

However, when the numerical approximation $\myvec{u}_h$, which is represented by $RT$ elements, is restricted to each triangular element $\mathcal{T}$, the components of $\myvec{u}_h|_\mathcal{T}$ are polynomials and thus differentiable \emph{within} $\mathcal{T}$. Hence, we consider the integral \eqref{eq:adv_integral} in an element-wise sense. Let $\Omega_h$ be a triangulation of the domain $\Omega\subset\mathbb{R}^N$. Then the integral involving the nonlinear advection term $\NL{u}{}$ can then be approximated by
\begin{equation}
\label{eq:mixed_NonlinearTerm}
\inprod{\NL{u}{h}}{\myvec{v}_h}_{\mathcal{T}} =
 \sum_{\mathcal{T}\in\Omega_h}\int_{\mathcal{T}}{\NL{u}{h}\cdot\myvec{v}_h}\intd{V},
 \quad\forall\myvec{v}_h\in RT_r(\Omega_h),
\end{equation}
where $RT_r(\Omega_h)$ is the Raviart-Thomas finite element space that
approximates $H(\tdiv;\Omega_h)$,
with $\nabla\myvec{u}_h$ defined in each element $\mathcal{T}\in\Omega_h$.
% --------------------------------------------------------------------------- 80
%
% --------------------------------------------------------------------------- 80
\subsection{Numerical Results on Manufactured Solutions}
\label{ssec:NonlinearAdvection_results}
% --------------------------------------------------------------------------- 80
In this subsection, we present convergence results for the full scheme:
PPE reformulation with the nonlinear advection term treated as described in
the previous subsection. To test the treatment of the advection term without
the complications from the PPE reformulation, we first study the spatial convergence
of the mixed formulation applied to the vector advection-diffusion
equation with electric boundary conditions. We then present the convergence
results for the full problem with the nonlinear term.

\subsubsection{Vector Nonlinear Advection-diffusion Equation with EBC}
To study the performance of the proposed treatment \eqref{eq:mixed_NonlinearTerm} of the advection term, we
consider the same semi-periodic domain and manufactured solution as in \Srm\ref{sssec:vhe_results}, but for the
vector-valued nonlinear advection-diffusion equation:
\begin{alignat}{2} \label{Eq:VNADE1}
\myvec{u}_t  + \nl{\myvec{u}}{}{} &= \Delta\myvec{u} + \myvec{f} \quad &&\mbox{for~}(x,y)\in(0,1)^2,\\ \label{Eq:VNADE2}
\myvec{n}\times\myvec{u} &= 0,\;\;\DIV\myvec{u}=0 \quad &&\mbox{for~}(x,y)\in[0,1)\times\{0,1\},\\ \label{Eq:VNADE3}
\myvec{u}(0,y) &= \myvec{u}(1,y) \quad &&\mbox{for~}0<y<1.
\end{alignat}

In the presence of the nonlinear advection term, systematic degradations in the spatial convergence order are observed (see Table~\ref{tab:spatial_convergence_vhe_nonlinear}). Specifically, non-convergent result is observed for $r=1$. For $r = 2,4$, the convergence rates are $r$, $r$, $r-1$, $r$ and $r-1$ for $\myvec{u}$, $\DIV\myvec{u}$, $\nabla\myvec{u}$, $\sigma$ and $\CURL\sigma$, respectively (in
both the $L^2$ and $L^\infty$ norms). The convergence rates are less clean when $r$ is odd, but follow a similar pattern as in the even order case.
These results indicate that degradation in the convergence order for the full PPE reformulation (\ref{eq:SR-PPE_momentum}--\ref{eq:SR-PPE_pressure}) should be expected.

\begin{table}[htb!]
	{\tiny
	\begin{tabular}{ccccccccccc}
	\toprule
		\multicolumn{11}{l}{Spatial approximation order $r = 1$, $(\sigma_h,\myvec{u}_h,p_h)\in P_1\times RT_0\times P_1$}\\
	\midrule
	$\dx$    & $\myvec{u}$ &  Rate  & $\DIV\myvec{u}$ & Rate & $\nabla\myvec{u}$ &  Rate  &  $\sigma$ & Rate &  $\CURL\sigma$ & Rate \\ 
	\midrule
	2.50E-01 & 6.93E+00 & ---- & 1.00E+01 & ---- & 1.37E+02 & ---- & 1.37E+02 & ---- & 2.76E+03 & ---- \\
	1.25E-01 & 3.28E+00 & 1.08 & 5.92E+00 & 0.76 & 1.37E+02 & 0.00 & 3.19E+01 & 2.10 & 1.45E+03 & 0.92 \\
	6.25E-02 & 1.67E+00 & 0.97 & 2.36E+00 & 1.33 & 1.37E+02 & 0.00 & 7.60E+00 & 2.07 & 7.25E+02 & 1.00 \\
	3.13E-02 & 8.75E-01 & 0.94 & 2.24E+00 & 0.07 & 1.37E+02 & 0.00 & 2.82E+00 & 1.43 & 3.65E+02 & 0.99 \\
	1.56E-02 & 5.07E-01 & 0.79 & 2.27E+00 & -0.02 & 1.37E+02 & -0.00 & 2.30E+00 & 0.30 & 1.87E+02 & 0.97\\
	7.81E-03 & 3.60E-01 & 0.49 & 2.28E+00 & -0.01 & 1.37E+02 & -0.00 & 2.27E+00 & 0.02 & 1.01E+02 & 0.88\\
	\midrule
	\multicolumn{11}{l}{Spatial approximation order $r = 2$, $(\sigma_h,\myvec{u}_h,p_h)\in P_2\times RT_1\times P_2$}\\
	\midrule
	$\dx$    & $\myvec{u}$ &  Rate  & $\DIV\myvec{u}$ & Rate & $\nabla\myvec{u}$ &  Rate  &  $\sigma$ & Rate &  $\CURL\sigma$ & Rate \\ 
	\midrule
	2.50E-01 & 2.13E+00 & ---- & 9.11E+00 & ---- & 9.33E+01 & ---- & 2.47E+01 & ---- & 1.25E+03 & ---- \\
	1.25E-01 & 7.92E-01 & 1.43 & 1.77E+00 & 2.37 & 6.26E+01 & 0.58 & 7.62E+00 & 1.70 & 4.03E+02 & 1.64 \\
	6.25E-02 & 2.14E-01 & 1.89 & 4.20E-01 & 2.07 & 3.28E+01 & 0.93 & 2.07E+00 & 1.88 & 1.17E+02 & 1.79 \\
	3.13E-02 & 5.48E-02 & 1.97 & 8.91E-02 & 2.24 & 1.66E+01 & 0.98 & 5.23E-01 & 1.98 & 3.92E+01 & 1.57 \\
	1.56E-02 & 1.38E-02 & 1.99 & 2.10E-02 & 2.08 & 8.33E+00 & 0.99 & 1.31E-01 & 2.00 & 1.63E+01 & 1.27 \\
	7.81E-03 & 3.45E-03 & 2.00 & 5.18E-03 & 2.02 & 4.17E+00 & 1.00 & 3.27E-02 & 2.00 & 7.65E+00 & 1.09 \\
	\midrule
	\multicolumn{11}{l}{Spatial approximation order $r = 3$, $(\sigma_h,\myvec{u}_h,p_h)\in P_3\times RT_2\times P_3$}\\
	\midrule
	$\dx$    & $\myvec{u}$ &  Rate  & $\DIV\myvec{u}$ & Rate & $\nabla\myvec{u}$ &  Rate  &  $\sigma$ & Rate &  $\CURL\sigma$ & Rate \\ 
	\midrule
	2.50E-01 & 1.33E+00 & ---- & 3.90E+00 & ---- & 7.10E+01 & ---- & 1.90E+01 & ---- & 6.41E+02 & ---- \\
	1.25E-01 & 1.28E-01 & 3.37 & 3.47E-01 & 3.49 & 1.67E+01 & 2.09 & 1.32E+00 & 3.85 & 8.12E+01 & 2.98 \\
	6.25E-02 & 1.58E-02 & 3.02 & 4.27E-02 & 3.02 & 4.24E+00 & 1.98 & 1.00E-01 & 3.72 & 1.16E+01 & 2.80 \\
	3.13E-02 & 2.03E-03 & 2.96 & 6.82E-03 & 2.65 & 1.06E+00 & 2.00 & 1.13E-02 & 3.16 & 2.25E+00 & 2.37 \\
	1.56E-02 & 2.90E-04 & 2.81 & 1.43E-03 & 2.25 & 2.66E-01 & 2.00 & 2.06E-03 & 2.45 & 5.15E-01 & 2.13 \\
	\midrule
	\multicolumn{11}{l}{Spatial approximation order $r = 4$, $(\sigma_h,\myvec{u}_h,p_h)\in P_4\times RT_3\times P_4$}\\
	\midrule
	$\dx$    & $\myvec{u}$ &  Rate  & $\DIV\myvec{u}$ & Rate & $\nabla\myvec{u}$ &  Rate  &  $\sigma$ & Rate &  $\CURL\sigma$ & Rate \\ 
	\midrule
	2.50E-01 & 2.42E-01 & ---- & 1.72E+00 & ---- & 1.94E+01 & ---- & 2.71E+00 & ---- & 1.57E+02 & ----\\ 
	1.25E-01 & 1.74E-02 & 3.80 & 7.45E-02 & 4.53 & 3.03E+00 & 2.68 & 1.80E-01 & 3.91 & 1.16E+01 & 3.76 \\
	6.25E-02 & 1.20E-03 & 3.86 & 2.53E-03 & 4.88 & 4.04E-01 & 2.91 & 1.40E-02 & 3.68 & 1.09E+00 & 3.41 \\
	3.13E-02 & 7.74E-05 & 3.96 & 1.33E-04 & 4.25 & 5.14E-02 & 2.97 & 9.28E-04 & 3.92 & 1.16E-01 & 3.22 \\
	1.56E-02 & 4.88E-06 & 3.99 & 7.79E-06 & 4.09 & 6.46E-03 & 2.99 & 5.89E-05 & 3.98 & 1.38E-02 & 3.07 \\
	\midrule
	\multicolumn{11}{l}{Spatial approximation order $r = 5$, $(\sigma_h,\myvec{u}_h,p_h)\in P_5\times RT_4\times P_4$}\\
	\midrule
	$\dx$    & $\myvec{u}$ &  Rate  & $\DIV\myvec{u}$ & Rate & $\nabla\myvec{u}$ &  Rate  &  $\sigma$ & Rate &  $\CURL\sigma$ & Rate \\ 
	\midrule
	2.50E-01 & 8.59E-02 & ---- & 4.32E-01 & ---- & 6.56E+00 & ---- & 1.28E+00 & ---- & 4.61E+01 & ---- \\
	1.25E-01 & 2.17E-03 & 5.31 & 8.36E-03 & 5.69 & 4.81E-01 & 3.77 & 2.93E-02 & 5.45 & 1.82E+00 & 4.66 \\
	6.25E-02 & 6.37E-05 & 5.09 & 1.74E-04 & 5.59 & 3.14E-02 & 3.94 & 5.27E-04 & 5.80 & 7.55E-02 & 4.59 \\
	3.13E-02 & 1.99E-06 & 5.00 & 7.18E-06 & 4.60 & 1.98E-03 & 3.99 & 1.28E-05 & 5.36 & 4.10E-03 & 4.20 \\
	1.56E-02 & 6.80E-08 & 4.87 & 3.68E-07 & 4.29 & 1.24E-04 & 4.00 & 5.41E-07 & 4.57 & 2.46E-04 & 4.06\\	
	\bottomrule
	\end{tabular}
}
\vspace{0.5em}
\caption{
	Observed spatial convergence rates in the $L^2$ norm for the vector nonlinear advection-diffusion equation
	with electric boundary conditions (\ref{Eq:VNADE1}--\ref{Eq:VNADE3}). The spatial approximation orders are
	$r=1,2,3,4,5$.}
\label{tab:spatial_convergence_vhe_nonlinear}
\end{table}

\subsubsection{The PPE Reformulation (\ref{eq:SR-PPE_momentum}--\ref{eq:SR-PPE_pressure}) (including the Nonlinear Advection Term)}
Consider the same manufactured solution, defined on the unit square
$\Omega=[0,1]^2$, as in \Srm\ref{sssec:linear_fullscheme}. However, since
here we solve the full problem (\ref{eq:SR-PPE_momentum}--\ref{eq:SR-PPE_pressure}), the forcing function is given by
$\myvec{f} = \myvec{u}_t+(\myvec{u}\cdot\nabla)\myvec{u} -
 \nu\Delta\myvec{u} + \nabla p$. We then use: the 3rd order spatial
discretization $(\sigma_h,\myvec{u}_h,p_h)\in P_3\times RT_2\times P_3$,
the 3rd order \imex{}(4,4,3) with time step $\dt = 0.2\dx$, the final time
$T = 3$, and $\lambda=30$ (the same as in the previous test cases).

\begin{table}[htb!]
	{\tiny
		\begin{tabular}{ccccccccccccccc}
			\toprule
			\multicolumn{15}{l}{Error convergence in $L^2$ norm}\\
			\midrule
			$\dx$   & $\myvec{u}$ &  & $\DIV\myvec{u}$ & & $\nabla\myvec{u}$ & & $\sigma$ & & $\CURL\sigma$ & & $p$ & & $\nabla p$   &      \\ 
			\midrule
			2.50E-01& 4.28E-02 & --- & 1.02E-01 & --- &  2.92E+00 & --- &  3.37E-01 & --- &  1.47E+01 & --- & 2.60E-01 & ---& 5.70E+00 & ---  \\
			1.25E-01& 8.14E-03 & 2.40& 9.39E-03 & 3.45&  1.03E+00 & 1.51&  3.52E-02 & 3.26&  2.97E+00 & 2.30& 5.33E-02 & 2.28& 9.49E-01 & 2.59\\
			6.25E-02& 1.05E-03 & 2.96& 7.22E-04 & 3.70&  2.62E-01 & 1.97&  4.74E-03 & 2.89&  5.93E-01 & 2.33& 1.45E-02 & 1.87& 2.42E-01 & 1.97\\
			3.13E-02& 1.36E-04 & 2.95& 5.89E-05 & 3.62&  6.57E-02 & 2.00&  8.15E-04 & 2.54&  1.37E-01 & 2.12& 3.69E-03 & 1.98& 6.11E-02 & 1.99\\
			1.56E-02& 1.84E-05 & 2.88& 5.91E-06 & 3.32&  1.64E-02 & 2.00&  1.64E-04 & 2.31&  3.34E-02 & 2.03& 9.18E-04 & 2.01& 1.52E-02 & 2.01\\
			7.81E-03& 2.70E-06 & 2.77& 6.84E-07 & 3.11&  5.64E-03 & 1.54&  3.68E-05 & 2.16&  8.30E-03 & 2.01& 2.28E-04 & 2.01& 3.82E-03 & 1.99\\
			\midrule
			\multicolumn{15}{l}{Error convergence in $L^\infty$ norm}\\
			\midrule
			$\dx$   & $\myvec{u}$ &  & $\DIV\myvec{u}$ & & $\nabla\myvec{u}$ & & $\sigma$ & & $\CURL\sigma$ & & $p$ & & $\nabla p$   &      \\ 
			\midrule
			2.50E-01&  2.37E-01 & ---&   3.20E-01 & ---&    1.30E+01 & ---&   2.32E+00 & ---&   1.05E+02 & ---&  1.24E+00 & ---& 2.87E+01 & --- \\
			1.25E-01&  1.12E-01 & 1.09&  7.30E-02 & 2.13&   7.80E+00 & 0.74&  2.31E-01 & 3.33&  2.30E+01 & 2.19& 3.51E-01 & 1.83& 5.27E+00 & 2.45 \\
			6.25E-02&  1.77E-02 & 2.65&  5.69E-03 & 3.68&   2.12E+00 & 1.88&  3.58E-02 & 2.69&  3.97E+00 & 2.54& 9.41E-02 & 1.90& 1.89E+00 & 1.48 \\
			3.13E-02&  2.34E-03 & 2.92&  3.56E-04 & 4.00&   5.36E-01 & 1.98&  4.99E-03 & 2.85&  9.46E-01 & 2.07& 2.34E-02 & 2.01& 5.60E-01 & 1.76 \\
			1.56E-02&  2.97E-04 & 2.98&  2.98E-05 & 3.58&   1.34E-01 & 2.00&  7.17E-04 & 2.80&  2.38E-01 & 1.99& 5.78E-03 & 2.02& 1.47E-01 & 1.93 \\
			7.81E-03&  3.73E-05 & 2.99&  3.15E-06 & 3.24&   3.35E-02 & 2.00&  1.46E-04 & 2.29&  6.02E-02 & 1.98& 1.43E-03 & 2.01& 3.73E-02 & 1.97 \\ 
			\bottomrule
		\end{tabular}
	}
\vspace{0.5em}
\caption{Error convergence for the full Navier-Stokes problem (including
	the nonlinear advection term) in the $L^2$ norm (top) and the $L^{\infty}$
	norm (bottom). The problem is solved with the 3rd order time-stepping
	\imex{}(4,4,3), the 3rd order spatial discretization
	$P_3\times RT_2\times P_3$, the final time $T = 3$, and $\lambda=30$.}
\label{tab:err_cvg_nonlinear}
\end{table}

The error convergence are shown in Table~\ref{tab:err_cvg_nonlinear}
in the $L^2$ norm (top) and in the $L^{\infty}$ norm (bottom). 
Convergent results are observed for all quantities but with slightly degraded convergence rates in comparison to the linear case in Table~\ref{tab:err_cvg_linear}. In particular, the convergence rate for $\DIV\myvec{u}$ is a little bigger than 3 as opposed to 4 in the linear case, and
$\CURL\sigma$, $p$ and $\nabla p$ show a 2nd order convergence.

% --------------------------------------------------------------------------- 80
%
% --------------------------------------------------------------------------- 80
\subsection{Numerical Results on Benchmark Tests}
\label{ssec:NonlinearAdvection_benchmark}
% --------------------------------------------------------------------------- 80
In this subsection, we demonstrate the performance of the proposed method
for solving two benchmark problems: lid-driven cavity in
\Srm\ref{sssec:NonlinearAdvection_benchmant_ldc}, and flow over a
backward-facing step in \Srm\ref{sssec:NonlinearAdvection_benchmant_bfs}.
The results show good agreement with reference data.
% --------------------------------------------------------------------------- 80
%
% --------------------------------------------------------------------------- 80
\subsubsection{Lid-driven Cavity}
\label{sssec:NonlinearAdvection_benchmant_ldc}
% --------------------------------------------------------------------------- 80
For the lid-driven cavity, we compute the flow in the unit square domain $[0,1]^2$, with
$\lambda = 10\/$, using the 3rd order spatial approximation
$P_3\times RT_2\times P_3$ and 3rd order \imex{} RK time-stepping (time step
$\dt = 0.8\dx$), on a regular triangular mesh with 16384 elements and a mesh
size $\dx = 1.5625\times 10^{-2}$. The velocity field is advanced forward in
time until it reaches steady state, for Reynolds numbers 100, 400 and 1000.
The flow starts at rest, with boundary conditions $\myvec{g}=(1,0)^T$ at
the top wall and no-slip elsewhere.

\begin{figure}[htb!]
\includegraphics[width=0.4\textwidth]{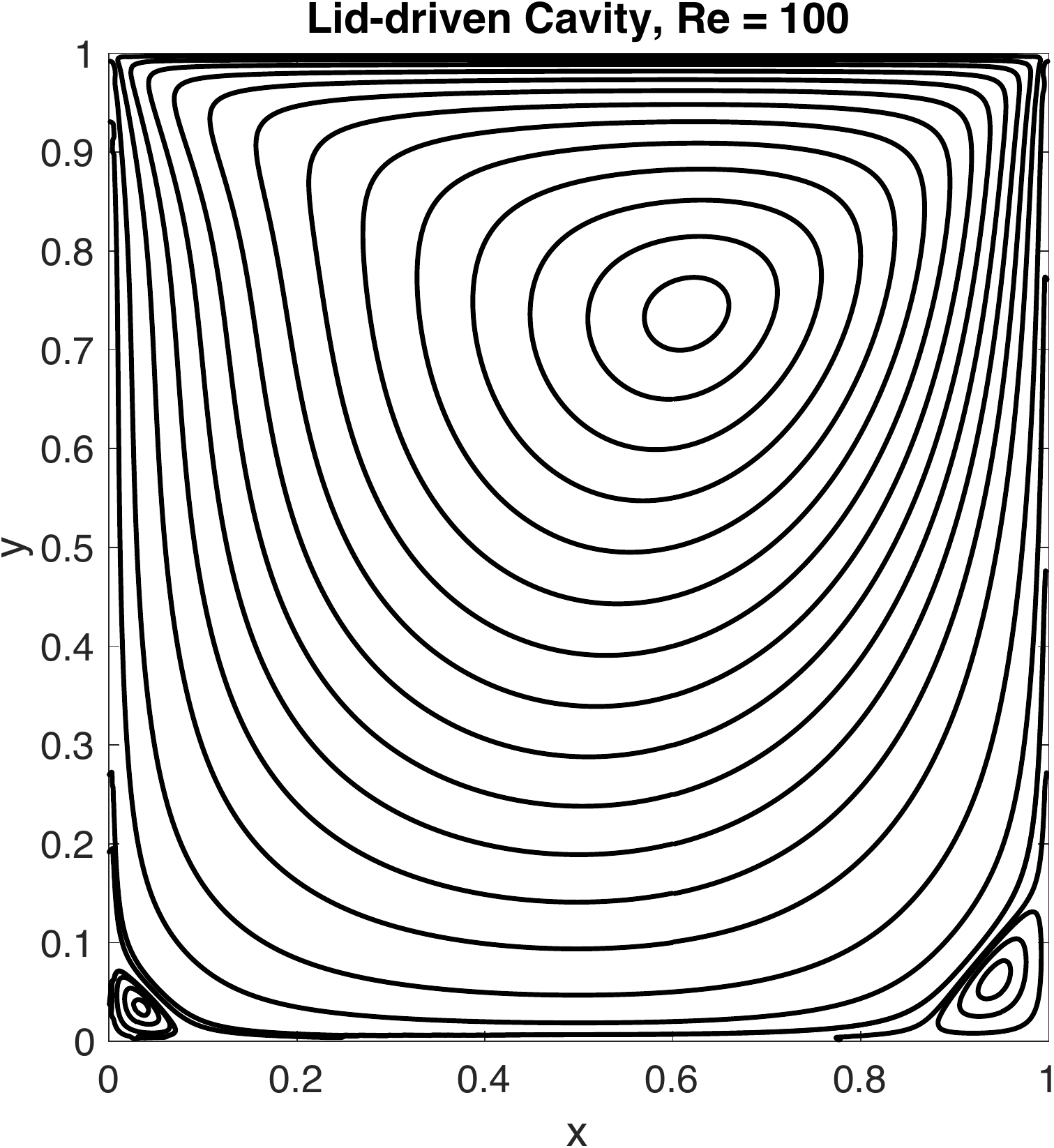}\qquad
\includegraphics[width=0.4\textwidth]{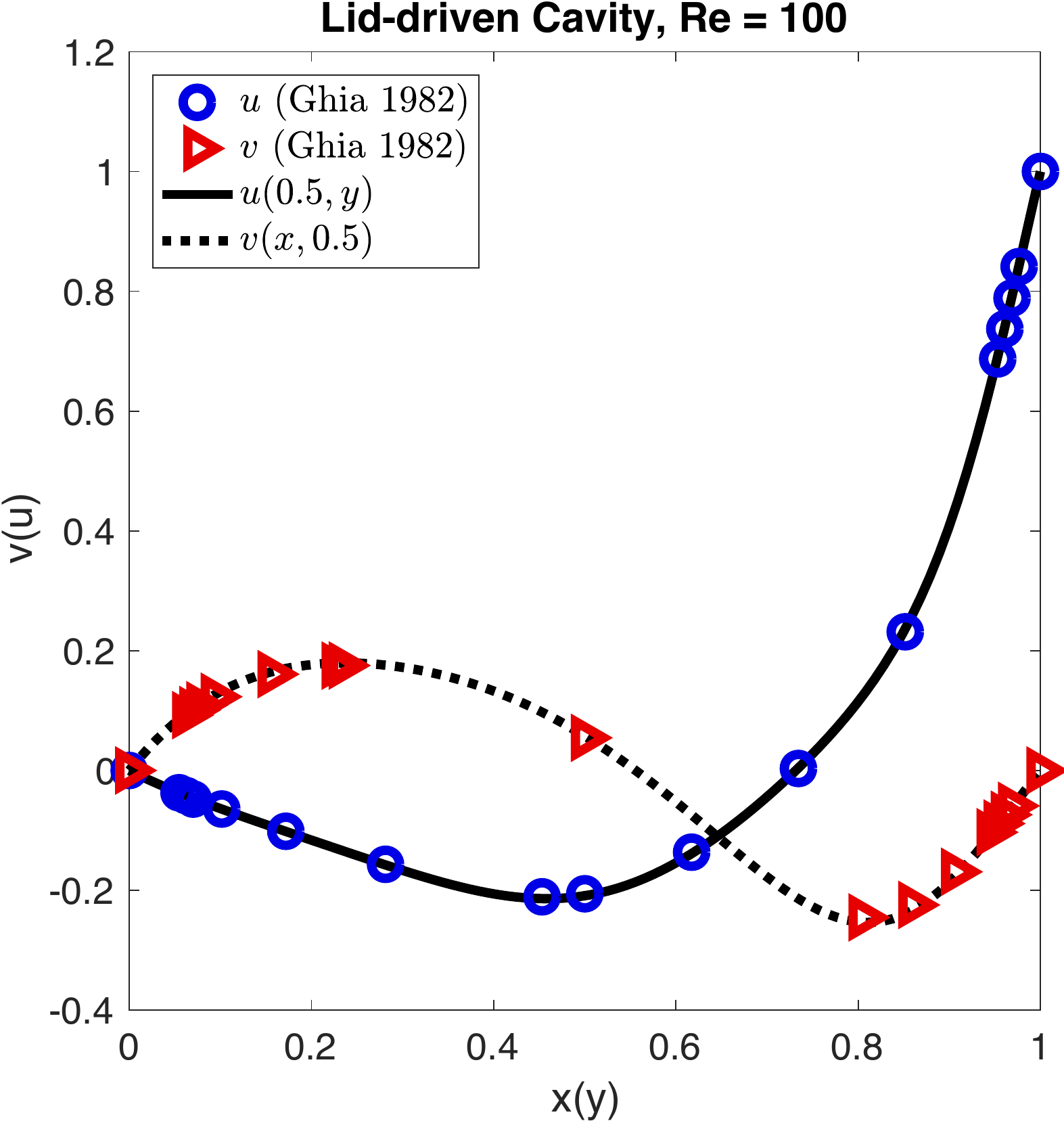}
\caption{Lid-driven cavity flow with $Re=100$.
Left: Streamlines at steady state.
Right: Velocity profiles along the centerlines (solid line: $u(0.5\/,\,y)$,
and dashed line: $v(x\/,\,0.5)$) compared with the reference data (blue
circles and red triangles) in \cite{GhiaGhiaShin1982}.}
\label{fig:ldc100}
\end{figure}

\begin{figure}[htb!]
\includegraphics[width=0.4\textwidth]{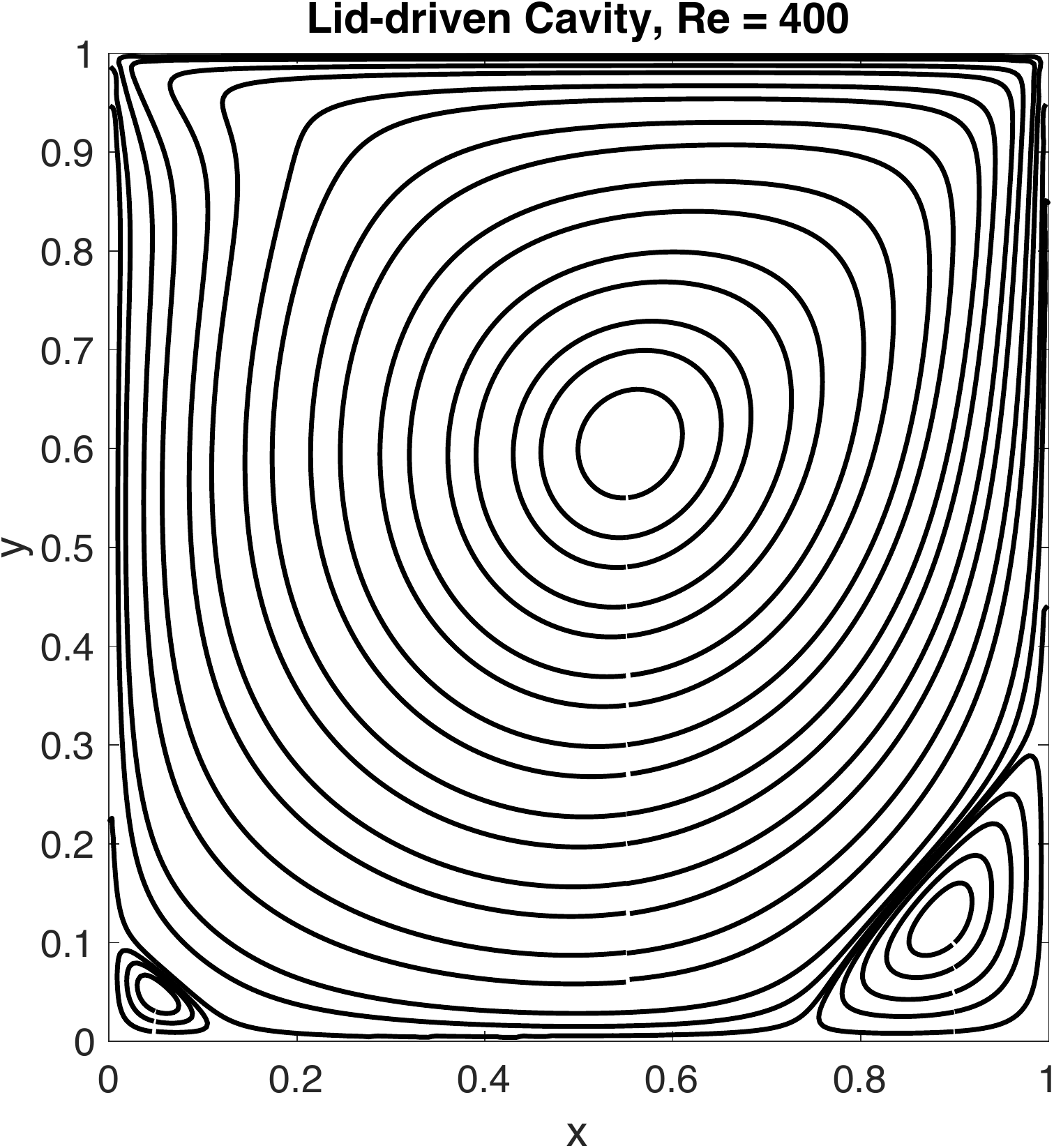}\qquad
\includegraphics[width=0.4\textwidth]{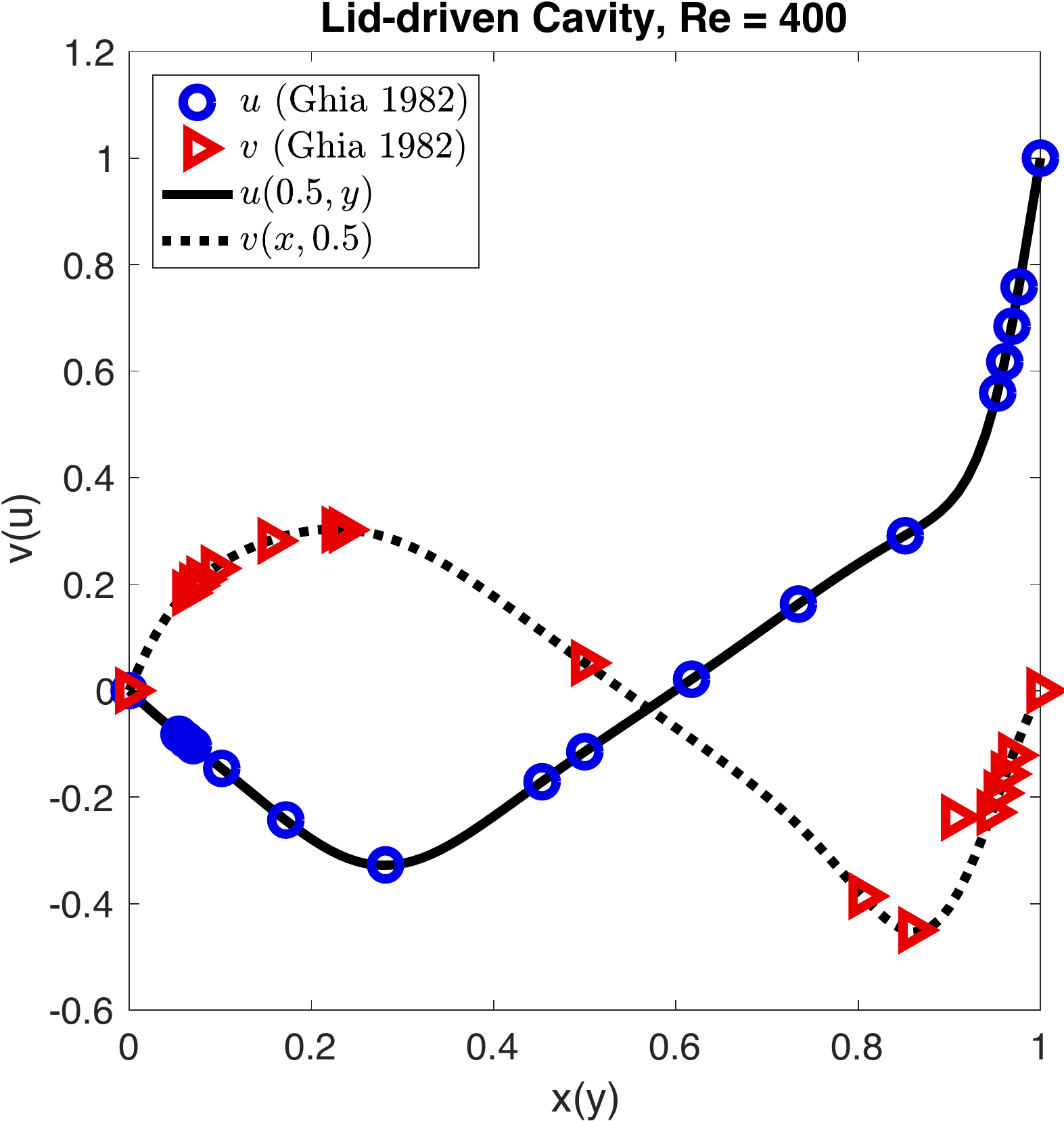}
\caption{Lid-driven cavity flow with $Re=400$.
Left: Streamlines at steady state.
Right: Velocity profiles along the centerlines (solid line: $u(0.5\/,\,y)$,
and dashed line: $v(x\/,\,0.5)$) compared with the reference data (blue
circles and red triangles) in \cite{GhiaGhiaShin1982}.}
\label{fig:ldc400}
\end{figure}

\begin{figure}[htb!]
\includegraphics[width=0.4\textwidth]{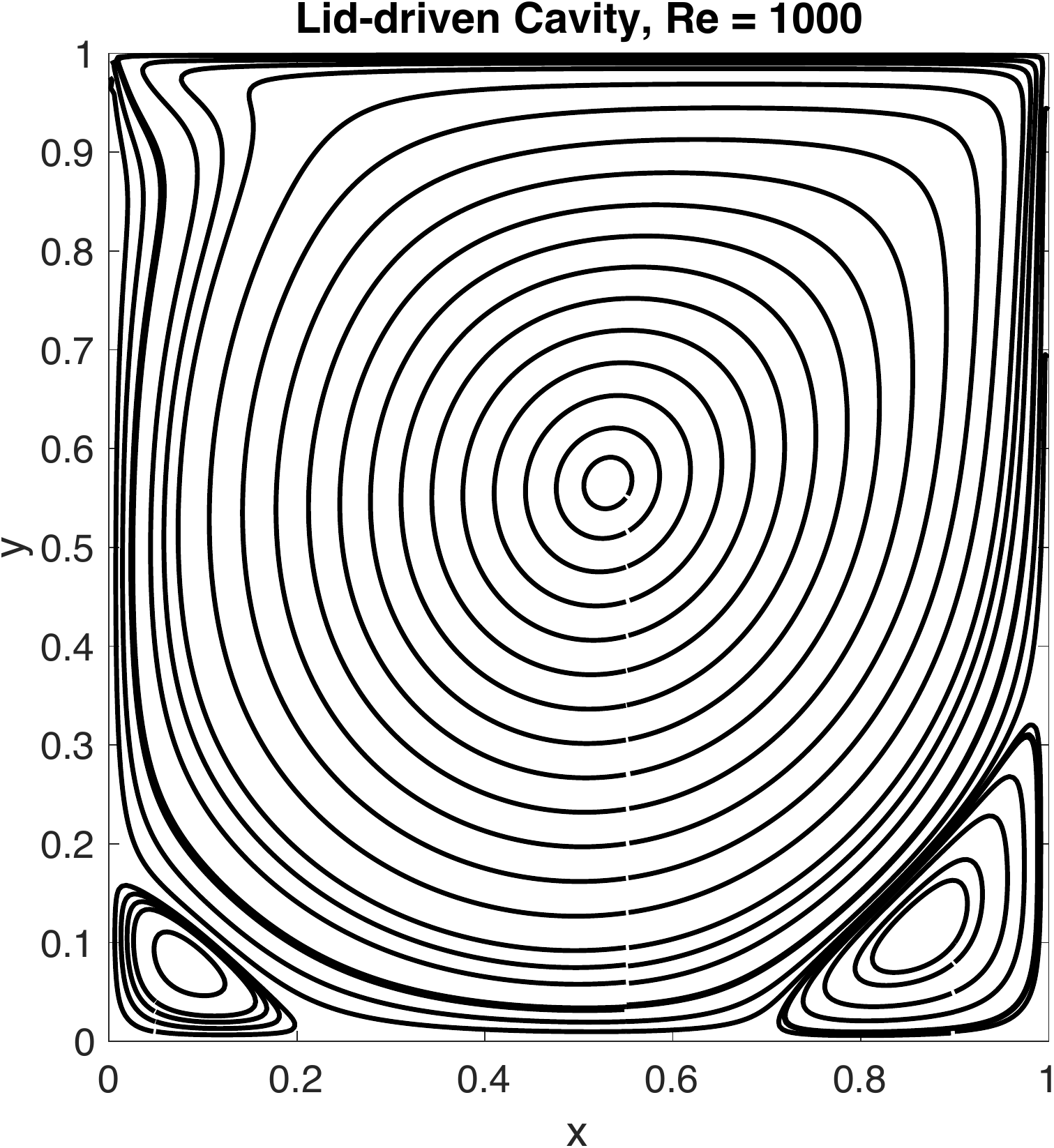}\qquad
\includegraphics[width=0.4\textwidth]{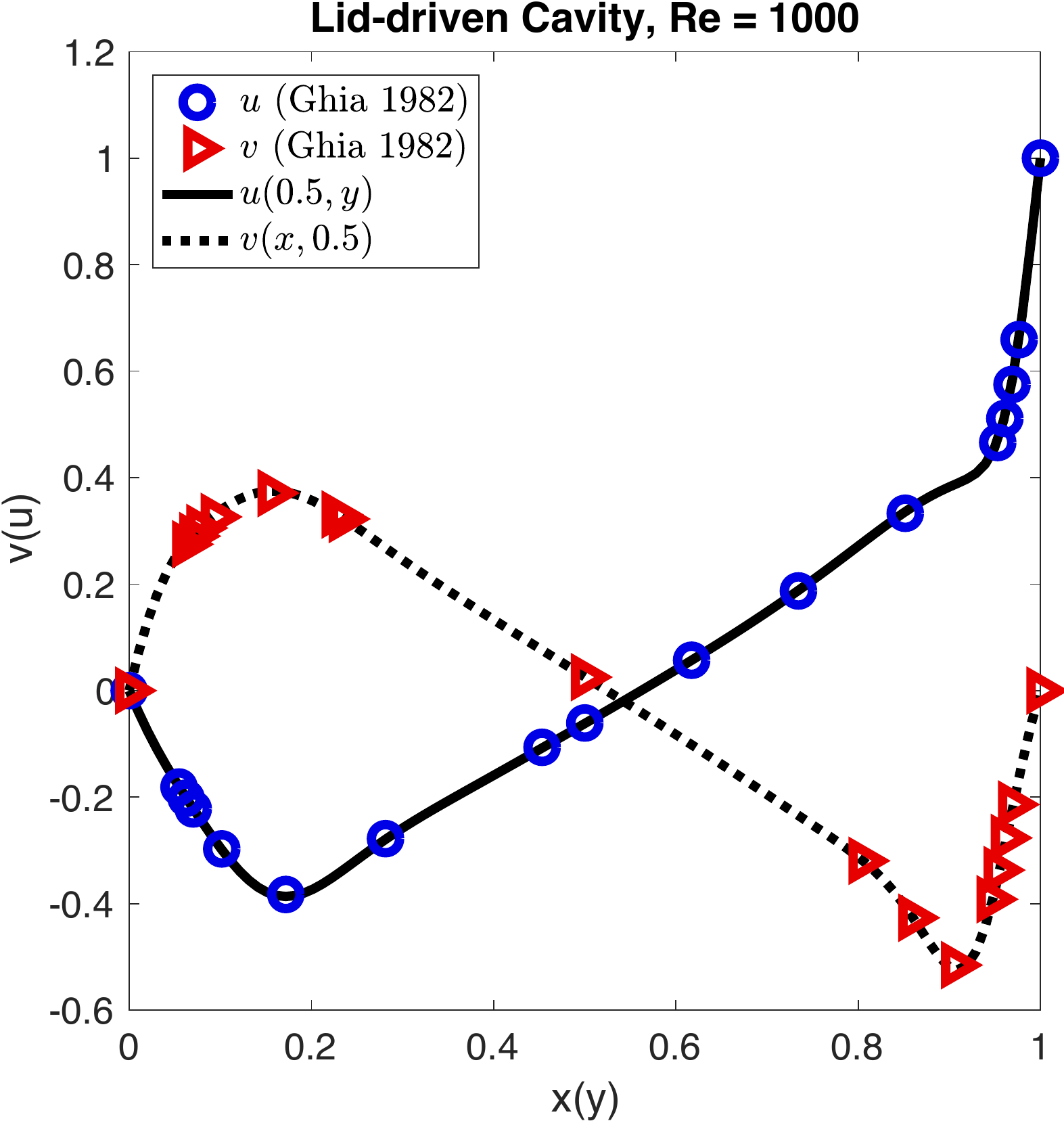}
\caption{Lid-driven cavity flow with $Re=1000$.
Left: Streamlines at steady state.
Right: Velocity profiles along the centerlines (solid line: $u(0.5\/,\,y)$,
and dashed line: $v(x\/,\,0.5)$) compared with the reference data (blue
circles and red triangles) in \cite{GhiaGhiaShin1982}.}
\label{fig:ldc1000}
\end{figure}

The results for $Re=100$, $400$ and $1000$ are shown in
Figures~\ref{fig:ldc100}--\ref{fig:ldc1000}. The streamlines for the steady
state flow are shown in the left panel of each figure.\footnote{\,The
   streamlines at the final time, $T$, are computed by numerically solving
   the ODE $\dfrac{\text{d}\myvec{x}}{\text{d}s}= \myvec{u}(\myvec{x},T)$
   with an explicit 4th order Runge-Kutta (RK4) scheme).}
The velocity profiles along the centerlines of the cavity (i.e.\,
$u(0.5\/,\,y)\/$ and $v(x\/,\,0.5)\/$) are plotted in the right panel,
together with the reference data from \cite{GhiaGhiaShin1982}.
Note that some streamlines in the plots end at the domain boundary, without
forming closed curves. This is due to the following facts: (i) the numerical
solution is not exactly divergence-free, (ii) the normal velocity condition
at the boundary is not enforced strongly (in the Dirichlet sense) but rather
through the ODE \eqref{eq:SR-normal}, which results in a small (as small as
the resolution) flow through the boundary.
The flow through the boundary is more pronounced at the top corners, where
discontinuities in the velocity occurs. Increasing $\lambda$ makes the
enforcing of the normal velocity condition stronger, and reduces the flow
through the boundary. However, making $\lambda\/$ too big would impose an undesirable
time step restriction, $\dt<O(\frac{1}{\lambda})$, through the relaxation term in \eqref{eq:SR-normal}.

% --------------------------------------------------------------------------- 80
%
% --------------------------------------------------------------------------- 80
\subsubsection{Backward-facing Step}
\label{sssec:NonlinearAdvection_benchmant_bfs}
% --------------------------------------------------------------------------- 80
In this test case, we compute the benchmark problem of flow 
over a backward-facing step for
$Re=100$ and $Re=200$. Again, we use the 3rd order spatial discretization
$P_3\times RT_2\times P_3$, a 3rd order \imex{} RK scheme and $\lambda=10$.
The computation uses a non-uniform triangular mesh with extra mesh refinement near the reentrant corner and the region behind the step. The minimum mesh size is $\dx = 1.7028\times 10^{-2}$ and the time step is set to $\dt = 0.02\dx$. In this case the domain of computation is
\begin{equation*}
	\Omega = [0,L]\times [-0.5,0.5]\backslash [0,0.5]\times[-0.5,0]\/,
\end{equation*}
where $L\/$ is the channel length, set to $L = 8\/$. No-slip boundary
conditions are imposed everywhere, except for the inflow and outflow
boundaries at $x=0$ and $x=L$. The inflow and outflow boundary conditions are
\begin{align*}
 \myvec{g}_{\text{inflow}} &= f(t) \left(12y(1-2y),0\right)^T,\\
 \myvec{g}_{\text{outflow}} &= f(t)\left(-3y^2+\dfrac{3}{4}, 0\right)^T\/,
\end{align*}
where $f(t)=1-e^{-6t^2}$, so that the flow is initially at rest and the
inflow and outflow increase gradually with time. The mean inflow velocity
$U$ reaches 1 for large enough $t$. We use the channel height $H=1$ as
the characteristic length, which gives a Reynolds number
$Re = HU/\nu = 1/\nu$.

\begin{figure}[htb!]
\centering
\includegraphics[width=0.8\textwidth]{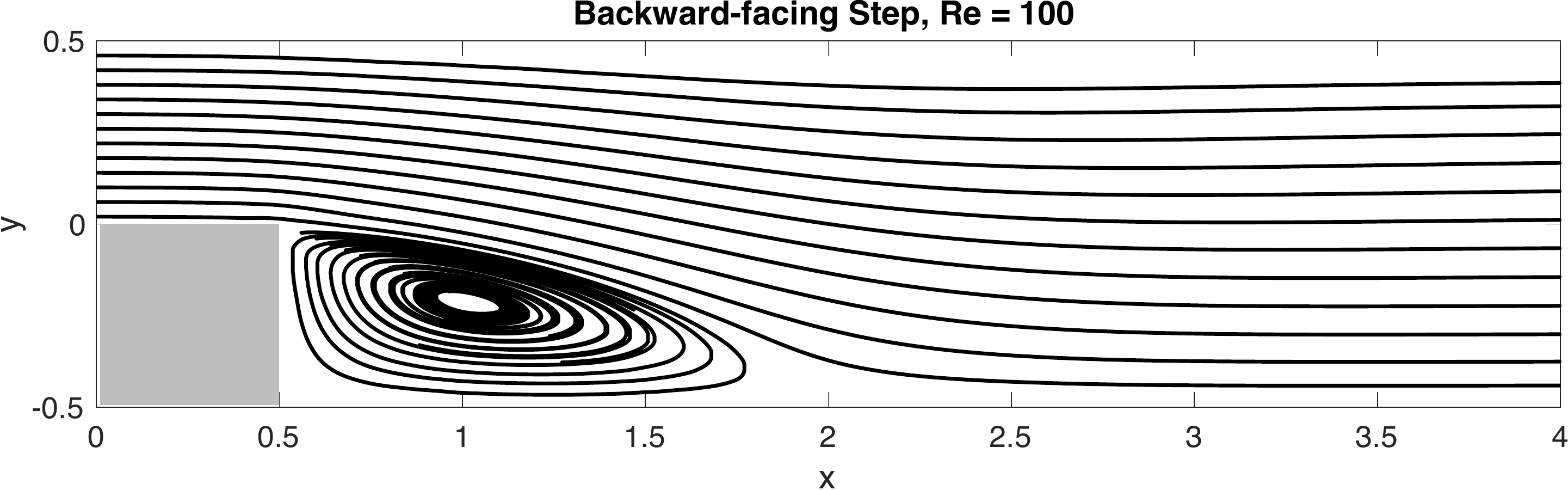}\\[0.5em]
\includegraphics[width=0.8\textwidth]{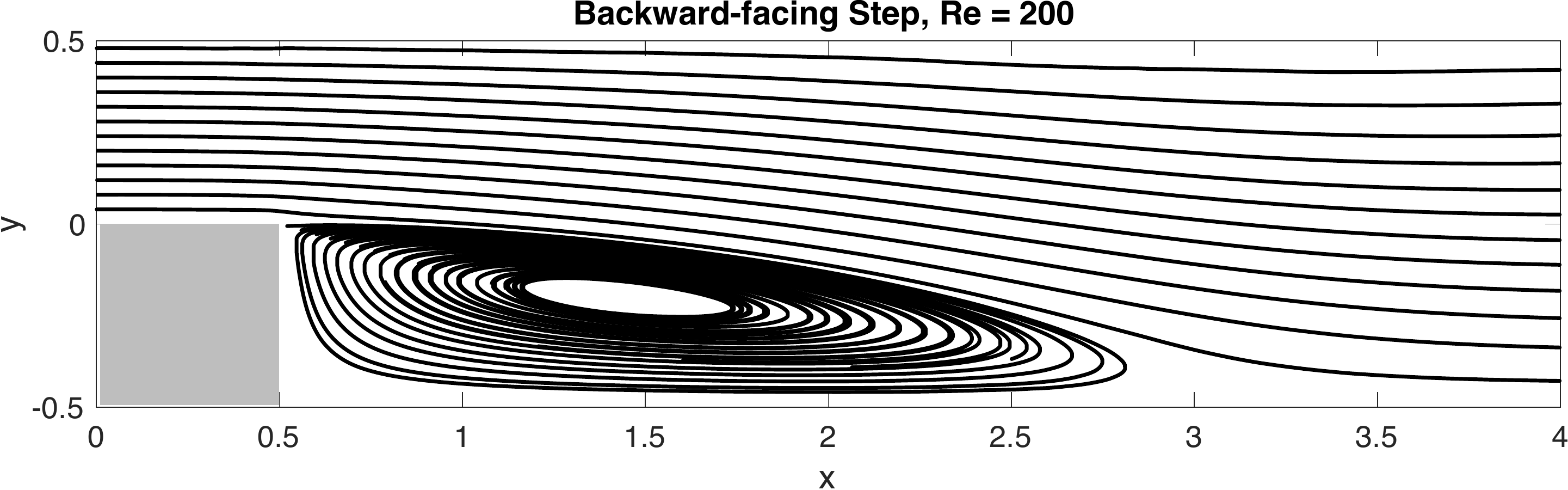}\\
\caption{Streamlines for the computation of a flow over a backward-facing
   step, at steady state, for $Re=100$ (top) and $Re=200$ (bottom).}
\label{fig:bfsstreamlines}
\end{figure}

For $Re=100$ and $200$, the only recirculating flow forms behind the step.
More regions of recirculating flow appear down the channel as the Reynolds
number increases. The streamlines shown in Figure~\ref{fig:bfsstreamlines}
were computed using the same procedure described in
\S\ref{sssec:NonlinearAdvection_benchmant_ldc}.

To compare our results with the reference data in \cite{Erturk2008} we use
the position of the reattachment point, where the line separating the
recirculating flow behind the step and the main flow in the channel meets
the channel wall, i.e.\, the domain boundary.
For this purpose introduce the nondimensional ratio $L_1/S$, where
$L_1$ is the distance between the foot of the step and the reattachment
point, and $S=0.5$ is the step height.
For $Re=100$, the ratio $L_1/S = 2.96$ in our computation compares well
with the reference data ratio: $2.922$.
For $Re=200$, our computation yields a ratio $L_1/S = 4.86$, while the
reference value is $4.982$. Our results show a rather good agreement
with the reference values.

% --------------------------------------------------------------------------- 80
%
% --------------------------------------------------------------------------- 80
\section{Conclusions and Outlook}
% --------------------------------------------------------------------------- 80
%
We investigated finite element formulations for a PPE reformulation of
the incompressible Navier-Stokes equations. In the PPE reformulation, the
momentum equation is in the form of a vector heat equation with electric
boundary conditions, and the pressure appears as a global function of the
velocity (obtained, at any time, as the solution to a Poisson equation). Thus
this reformulation allows for high-order time-stepping via standard schemes.
In particular, the decoupling of the velocity and the pressure can be
achieved by \imex{} time-stepping schemes, and modern high-order \imex{} schemes
can be applied in a straightforward manner.

We demonstrated via numerical examples that our proposed numerical schemes,
based on the PPE reformulation (\ref{eq:SR-PPE_momentum}--\ref{eq:SR-PPE_pressure}), have the potential to achieve high-order
both in space and in time, while avoiding severe time step restrictions. The
methods have important advantages:
\begin{enumerate}[(i)]
 \item The use of an \imex{} time-stepping strategy decouples the velocity and
 the pressure in the numerical methods, and at the same time avoids diffusive
 time step restrictions.
 \item The methods can achieve high-order in time with off-the-shelf
 high-order \imex{} RK schemes. Therefore the code can be easily adapted to
 newly developed \imex{} schemes with better properties, for instance,
 unconditional stability, order reduction avoidance, etc. A natural future
 research direction includes developing new \imex{} RK time-stepping schemes
 for the PPE reformulations with such properties.
 \item  Standard mixed finite element formulations of incompressible fluid
 flow problems (Stokes and Navier-Stokes equations) require the velocity
 and pressure approximations to satisfy the inf-sup condition for
 stability, which limits the choices of finite element approximations. Our
 new schemes allow for more flexible choices of finite element spaces for
 the velocity and pressure, avoiding the inf-sup condition.
\end{enumerate}
Note also that an interesting feature of the methods studied in this
paper is that the quantity $\nabla \cdot \myvec{u}$ converges to zero at an additional
order relative to the accuracy of the velocity field itself. As a
consequence, for well-resolved computations, the methodology will
generally yield velocity fields that are extremely close to
divergence-free, even though no discrete incompressibility principle
needs to be imposed.

Despite the important advantages mentioned above, the methods also
have some limitations:
\begin{enumerate}[(i)]
\item The discretization of the (nonlinear) advection term is ad-hoc, and
does not fit into the finite element formulation due to the discontinuities
across elements in the tangential velocity. The current approach leads to
convergent methods but results in degradations in the error convergence
rates.
\item The mixed formulation for the velocity deals with the EBC naturally.
However, it introduces a saddle point problem. Therefore the choices for
the FE approximations for $\myvec{u}$ and $\sigma$ need to satisfy the
inf-sup condition. The mixed formulation also increases the degrees of
freedom of the discrete problem, as the new variable
$\sigma = \CURL\myvec{u}$ is introduced.
\end{enumerate}

\appendix
% --------------------------------------------------------------------------- 80
\section{Discrete solvability of the Poisson equation}\label{Sec:AppendixPoisson}
% --------------------------------------------------------------------------- 80
The pressure Poisson problem \eqref{eq:linear_SR-PPE_p} has pure Neumann
boundary conditions prescribed, hence it is solvable only if a compatibility condition is satisfied (see Remark~\ref{rem:PPE_solvability}), and the solution is unique only up to an additive constant. In order to single out a unique solution, and at the same time to obtain a stable approximate solution in case the compatibility condition is not exactly satisfied (due to approximation errors), we employ the following standard least-squares approximation procedure. First, we impose an extra zero-mean constraint on the pressure:
\begin{equation*}
	\int_{\Omega} p \intd{V} = 0\/.
\end{equation*}
This gives rise to an augmented system for the pressure with an additional
scalar variable representing the Lagrange multiplier for the zero-mean
constraint. For simplicity, the presentation here is carried out for homogeneous boundary conditions
$\myvec{g}=0$. The weak formulation then becomes: Find
$p\in H^1(\Omega)$ and $c\in\mathbb{R}$ s.t.
\begin{subequations}
	\label{eq:augmented_PPE}
	\begin{alignat}{2}
		\inprod{\nabla p}{\nabla q} + \inprod{c}{q} &= \inprod{\myvec{f}}{\nabla q}
		- \nu\bdryint{\myvec{n}\cdot(\CURL\sigma)q} +\lambda\bdryint{(\myvec{n}\cdot\myvec{u})q}
		&&\quad\forall q\in H^1(\Omega),\\
		\inprod{p}{d} & =0 &&\quad\forall d\in\mathbb{R}.
	\end{alignat}
\end{subequations}
Using standard nodal-based finite elements yields a linear system of the form
\begin{equation}
	\begin{pmatrix}
		K & \myvec{r} \\ \myvec{r}^T & 0
	\end{pmatrix}
	\cdot
	\begin{pmatrix}
		\myvec{P}\\ c
	\end{pmatrix}
	=
	\begin{pmatrix}
		\myvec{F}\\ 0
	\end{pmatrix} .
	\label{eq:augmented_PPE_linsys}
\end{equation}
Here $K$ is the (symmetric) stiffness matrix of the FEM discretization of the Laplacian operator, $\myvec{r}$ is the vector corresponding to constant functions (the null-vector of $K$), $\myvec{P}$ is the solution vector for the pressure $p$, and $c$ is the Lagrange multiplier. While the stiffness matrix $K$ is singular, the augmented matrix in \eqref{eq:augmented_PPE_linsys} is nonsingular and it yields the following solution. Left-multiplying the equation $K\myvec{P}+c\myvec{r} = \myvec{F}$ by $\myvec{r}^T$ implies that $c = (\myvec{r}^T\myvec{F})/(\myvec{r}^T\myvec{r})$, thus one has $K\myvec{P} = \textrm{proj}_{\mathcal{R}(K)} \myvec{F}$, where $\mathcal{R}(K)$ is the range of $K$; and $\myvec{r}^T \myvec{P}=0$ restricts that $\myvec{P} \in \mathcal{R}(K)$.

Hence, the augmented system \eqref{eq:augmented_PPE} addresses the two issues arising in solving the pressure Poisson equation: (i)~it fixes the additive constant in $p$ by choosing the zero-mean solution, and (ii)~it ensures solvability, even when $\DIV\myvec{u}=0$ is violated (see Remark~\ref{rem:PPE_solvability}), by projecting the right hand side $\myvec{F}$ to the range of $K$.
Note that the augmented system \eqref{eq:augmented_PPE_linsys} is similar to the one discussed in \cite{Henshaw1994,HenshawPetersson2003,ShirokoffRosales2010} where finite difference approaches are used.

% --------------------------------------------------------------------------- 80
% --------------------------------------------------------------------------- 80
% --------------------------------------------------------------------------- 80
\section*{Acknowledgment}
% --------------------------------------------------------------------------- 80
The authors wish to acknowledge support by the National Science Foundation
through the grants
DMS--1719637 (Rosales),
DMS--1719640 (Seibold and Zhou), and
DMS--1719693 (Shirokoff).
In addition, the authors received partial support through the grants
DMS--1614043 (Rosales), and
DMS--1318709 (Seibold and Zhou).

% --------------------------------------------------------------------------- 80
%\appendix
% --------------------------------------------------------------------------- 80
% --------------------------------------------------------------------------- 80
% --------------------------------------------------------------------------- 80
% --------------------------------------------------------------------------- 80
% --------------------------------------------------------------------------- 80
%\newpage
% --------------------------------------------------------------------------- 80
% \newpage
\bibliographystyle{plain}
\bibliography{../ppe_references_complete}

\end{document}